\documentclass{amsart}

\usepackage[mathscr]{eucal}
\usepackage{amssymb}
\usepackage{upgreek}
\usepackage{graphicx}
\usepackage{psfrag} 
\usepackage[dvipsnames]{xcolor}
\usepackage[normalem]{ulem}
\usepackage{amsthm}
\usepackage{bbm}
\usepackage{enumerate}
\usepackage[color=yellow]{todonotes}
\usepackage{array}
\usepackage{amsmath}

\usepackage{hyperref}
\hypersetup{colorlinks=true,linkcolor={Brown},citecolor={Brown},urlcolor={Brown}}
\usepackage{mathtools}

\usepackage{cleveref}
\crefname{Thm}{Theorem}{Theorems}
\crefname{Rem}{Remark}{Remarks}
\crefname{Prop}{Proposition}{Propositions}
\crefname{Cor}{Corollary}{Corollaries}
\crefname{Cons}{Construction}{Constructions}
\crefname{Exa}{Example}{Examples}
\crefname{Lem}{Lemma}{Lemmas}
\crefname{Rec}{Recollection}{Recollections}
\crefname{subsection}{Subsection}{Subsections}

\numberwithin{equation}{section}
\makeindex

\usepackage[all]{xy}
\CompileMatrices

\newdir{ >}{{}*!/-10pt/\dir{>}}

\swapnumbers 

\newtheorem{Thm}[equation]{Theorem}
\newtheorem{Prop}[equation]{Proposition}
\newtheorem{Lem}[equation]{Lemma}
\newtheorem{Cor}[equation]{Corollary}

\theoremstyle{remark}
\newtheorem{Rem}[equation]{Remark}
\newtheorem{Def}[equation]{Definition}

\newtheorem{Not}[equation]{Notation}

\newtheorem{Cons}[equation]{Construction}

\newtheorem{Rec}[equation]{Recollection}

\theoremstyle{definition}
\newtheorem*{Ack}{Acknowledgements}
\newtheorem*{AI}{Statement about use of A.I}

\newcommand{\nc}{\newcommand}
\nc{\dmo}{\DeclareMathOperator}

\dmo{\End}{End}
\dmo{\Fun}{\mathsf{Fun}} 
\dmo{\Hom}{Hom}
\dmo{\incl}{incl}
\dmo{\Ind}{Ind}
\dmo{\ind}{ind}
\dmo{\modcat}{mod} 
\dmo{\Modcat}{Mod} 
\dmo{\Mor}{Mor}%
\dmo{\Obj}{Obj}

\nc{\gr}{\mathsf{Gr}}

\newcommand{\cI}{\mathcal I}

\newcommand{\cK}{\mathcal K}

\newcommand{\e}{\mathcal E}
\newcommand{\s}{\mathcal S}
\newcommand{\n}{\mathcal N}
\newcommand{\w}{\mathcal W}

\newcommand{\cC}{\mathcal C}

\NewDocumentCommand{\Alp}{o}{\IfNoValueTF{#1}{A}{#1}_{\scriptscriptstyle+}}

\NewDocumentCommand{\Bup}{o}{\IfNoValueTF{#1}{B}{#1}^{\scriptscriptstyle+}}

\dmo{\lp}{Lplus}
\dmo{\Up}{Uplus}

\dmo{\Mod}{Mod}

\dmo{\Res}{Res}
\dmo{\Inf}{Inf}
\dmo{\Iso}{Iso}
\dmo{\Co}{Coind}
\dmo{\Deff}{Def}

\dmo{\res}{res}
\dmo{\con}{con}
\dmo{\iso}{iso}
\dmo{\deff}{def}

\newcommand{\lind}[1]{{}_{#1}\!\ind}
\newcommand{\lres}[1]{{}_{#1}\!\res}
\newcommand{\ldef}[1]{{}_{#1}\!\deff}
\newcommand{\linf}[1]{{}_{#1}\!\inf}
\newcommand{\liso}[1]{{}_{#1}\!\iso}

\newcommand{\lInd}[1]{{}_{#1}\!\Ind}
\newcommand{\lRes}[1]{{}_{#1}\!\Res}

\newcommand{\lInf}[1]{{}_{#1}\!\Inf}

\newcommand{\lCo}[1]{{}_{#1}\!\Co}

\newcommand{\iInd}[2]{{}_{#1}\!\Ind_{#2}}
\newcommand{\iRes}[2]{{}_{#1}\!\Res_{#2}}
\newcommand{\iInf}[2]{{}_{#1}\!\Inf_{#2}}

\newcommand{\iCo}[2]{{}_{#1}\!\Co_{#2}}
\newcommand{\iDef}[2]{{}_{#1}\!\Deff_{#2}}

\newcommand{\ires}[2]{{}_{#1}\!\res_{#2}}

\renewcommand{\inf}{\mathop{\mathrm{inf}}\nolimits}

\dmo{\mob}{m\ddot{o}b}
\dmo{\mk}{mk}
\dmo{\km}{km}
\dmo{\die}{die}
\dmo{\can}{can}
\dmo{\lin}{lin}
\dmo{\tm}{tm}
\dmo{\Img}{Im}
\dmo{\Ker}{Ker}
\dmo{\ev}{ev}

\newcommand{\br}{\text{br}}

\dmo{\tr}{tr}

\dmo{\hm}{Hom}
\dmo{\en}{End}

\dmo{\id}{id}
\dmo{\op}{op}
\dmo{\Sd}{Sd}
\dmo{\ch}{ch}

\dmo{\SL}{SL}

\nc{\To}{\Rightarrow}
\nc{\xTo}[1]{\overset{#1}{\,\To\,}}
\nc{\isoTo}{\overset{\sim}{\To}}
\nc{\isoto}{\overset{\sim}{\to}}
\nc{\Too}{\Longrightarrow}
\nc{\onto}{\twoheadrightarrow}

\nc{\MS}[1]{(MS\,\ref{it:MS-#1})}
\nc{\GC}[1]{(GC\,\ref{it:GC-#1})}

\newcommand{\ol}[1]{\overline{#1}}

\nc{\congby}[1]{\overset{\textrm{#1}}{\cong}}

\nc{\grfun}[1]{\Fun^{\mathsf{Gr}}(#1)}

\nc{\eps}{\varepsilon}
\nc{\leps}{{}^{\ell}\eps}
\nc{\leta}{{}^{\ell}\eta}
\nc{\etalp}{{}_{+}\eta}
\nc{\epslp}{{}_{+}\eps}
\nc{\etaup}{{}^{+}\eta}
\nc{\epsup}{{}^{+}\eps}
\nc{\Hatice}[1]{{\color{RoyalBlue}#1}}

\title{Plus Constructions and Canonical Induction for Green Categories}
\author{Hatice Mutlu}
\thanks{\url{https://mutluhtc.github.io}}
\date{\today}


\begin{document}
\maketitle
\small

\begin{abstract}\normalsize
We develop canonical induction for Green categories associated with
Mackey systems. We establish a change-of-rings framework in which the
lower-plus and upper-plus constructions are realized by induction and
coinduction of modules. The two plus adjunctions yield the mark
morphism, and, when the relevant group orders are invertible in the
coefficient ring, M\"obius inversion gives an explicit inverse.
Using the resulting linearization and tom Dieck maps, we construct a
canonical induction map, give an explicit formula for it, and
characterize when it is a right inverse of linearization. In the 
classical setting of Mackey functors, our construction and splitting 
criterion recover Boltje's results.
\end{abstract}

\tableofcontents

\section{Introduction}
\label{sec:introduction}

Brauer's induction theorem states that every complex character of a
finite group can be expressed as an integral linear combination of
characters induced from linear characters of elementary subgroups.
Boltje \cite{Boltje98a} developed canonical induction in the general
setting of Mackey functors. In this framework, a Mackey functor $M$
plays the role of the character ring, while the restriction functor
$A$ plays the role of the linear characters. Unlike $M$, the functor
$A$ is equipped with restriction maps but not induction maps.
The lower-plus construction produces a Mackey functor $A_+$ from
$A$, together with an induction morphism $A_+\to M$.
Canonical induction seeks a restriction-compatible section
to this morphism, and Boltje established criteria for when
such a splitting exists.
 The upper-plus construction provides a complementary construction
$A^+$, related to $A_+$ by the mark morphism $A_+\to A^+$.
This relationship plays a central role in constructing canonical
induction formulas. 

Co\c{s}kun \cite{Coskun07} later gave a
module-theoretic interpretation of the plus constructions in terms
of induction and coinduction. Co\c{s}kun subsequently adapted these 
module-theoretic methods to biset functors on the subquotients of a 
fixed finite group, studying induction, coinduction, and 
a mark morphism connecting these constructions \cite{Coskun08}.
Boltje, Raggi-C\'ardenas and Valero-Elizondo
\cite{BoltjeRaggiCardenasValeroElizondo19}
extended the lower- and upper-plus constructions to biset
functors on suitable subcategories of the biset category.
They also established the corresponding adjunctions and
studied the mark morphism and the multiplicative structures
of these constructions for Green biset functors. Further
extensions include the fibered biset-functor constructions of
Calder\'on \cite{Calderon25} and the categorification using Mackey
$2$-functors developed by Balmer and the author
\cite{BalmerMutlu57pp}.

The purpose of the present paper is to study the plus constructions
and canonical induction in the setting of Green categories associated
with Mackey systems. Let $\mathfrak K$
be a set of finite groups closed under taking subgroups.
A Mackey system on $\mathfrak K$ (\cref{rec:mackey-systems}) is a subcategory of the category of
finite groups with object set $\mathfrak K$, whose morphisms satisfy
the Mackey-system axioms. Let $(\mathcal I,\mathcal R)$ be a pair of
Mackey systems satisfying the Green-category axioms in \cref{rec:green-categories}.
Following Barker \cite{Barker16b}, a Green category on $\mathfrak K$ associated with 
$(\mathcal I,\mathcal R)$ is defined as the
$\mathbb Z$-linear subcategory of the biset category generated by
induction along morphisms of $\mathcal I$ and restriction along
morphisms of $\mathcal R$.

We concentrate on ordinary induction, namely induction along
injective homomorphisms, while allowing restriction along the
morphisms of an arbitrary Mackey system.

\begin{Thm}(\Cref{sec:green-categories-from-mackey-systems})
\label{thm:intro-green-category}
Let $\mathcal R$ be a Mackey system on $\mathfrak K$, and let
$\mathcal I$ be its subcategory of injective morphisms. Then
$(\mathcal I,\mathcal R)$ satisfies the Green-category axioms
and hence defines a Green category.

Conversely, if a pair $(\mathcal I,\mathcal R)$ of Mackey systems
defines a Green category and every morphism of $\mathcal I$
is injective, then $\mathcal I$ consists precisely of the
injective morphisms of $\mathcal R$.
\end{Thm}

Thus the induction system and the associated Green category are
determined by the restriction Mackey system $\mathcal R$.
We compare this setting with the admissible-biset framework of
\cite{BoltjeRaggiCardenasValeroElizondo19} in
\Cref{rem:admissible-biset-categories}.

Recall that a group category is a $\mathbb Z$-linear subcategory
of the biset category whose morphism groups have bases consisting
of transitive bisets. Let $\mathcal R$ be a Mackey system on $\mathfrak K$, and let
$\mathcal I$ be its subcategory of injective morphisms. 
Associated with $(\mathcal I, \mathcal R)$, we define four group
categories $\mathcal S,\mathcal W,\mathcal E,\mathcal N$ 
(\cref{ssec:four-group-categories}), with
inclusions 
\[
\xymatrix@C=1ex@R=1ex{
	& \mathcal N & \\
	\mathcal W \ar@{}[ur]|-{\rotatebox{45}{$\subseteq$}}
	&& \mathcal E \ar@{}[ul]|-{\rotatebox{135}{$\subseteq$}}\\
	& \mathcal S  \ar@{}[ul]|-{\rotatebox{135}{$\subseteq$}} \ar@{}[ur]|-{\rotatebox{45}{$\subseteq$}} &
}
\]
Here $\mathcal N$ is the Green category associated with
$\mathcal R$,  the category $\mathcal E$ is generated by its restriction morphisms,
the category $\mathcal W$ is generated by ordinary induction and restriction
along surjective homomorphisms, and the category $\mathcal S$ is generated by
the surjective-restriction morphisms. The category
$\mathcal N$ is the inflation Mackey category of
\cite{Barker16a}. We identify the group categories with their locally unital
quiver algebras, as explained in \Cref{ssec:group-categories}.

Let $R$ be a commutative unital ring.
For each group category $\mathcal X$, we write
$R\mathcal X=R\otimes_{\mathbb Z}\mathcal X$
for its extension of scalars to $R$.
As explained in \Cref{ssec:group-categories}, we use the same
notation for the resulting $R$-linear category and its locally
unital quiver algebra, and identify functors on these categories
with locally unital modules over the corresponding algebras.
In this setting, we obtain the following analogue
of Co\c{s}kun's Mackey structure theorem
\cite[Thm.~3.2]{Coskun07}.
\begin{Prop}(\Cref{sec:four-group-categories})
\label{prop:intro-bimodule-decomposition}
Composition induces an isomorphism of
$(R\mathcal W,R\mathcal E)$-bimodules
\[
R\mathcal W\otimes_{R\mathcal S}R\mathcal E
\xrightarrow{\ \sim\ } R\mathcal N.
\]
\end{Prop}
We also consider Green $R\mathcal S$-,
$R\mathcal E$-, and $R\mathcal N$-functors.
Their values are unital $R$-algebras, and their
restriction maps are unital algebra homomorphisms.
Green $R\mathcal N$-functors additionally satisfy
the Frobenius relations for induction, as in
\cite{Barker16b}.

Following the module-theoretic approach of Co\c{s}kun
\cite{Coskun07}, we define the lower-plus construction
by induction from $R\mathcal E$ to $R\mathcal N$ 
(\cref{con:lower-plus})
and the upper-plus construction by coinduction from
$R\mathcal W$ to $R\mathcal N$ after inflation from
$R\mathcal S$ to $R\mathcal W$ (\cref{con:upper-plus}).
In \cref{prop:lower-plus-cofixed} and \cref{prop:upper-plus-fixed}, 
we identify these constructions with the classical
cofixed and fixed descriptions of Boltje \cite{Boltje98a} .
\begin{Thm}(\Cref{sec:plus-constructions})
\label{thm:intro-plus-constructions}
The lower-plus construction defines a functor
\[
(-)_+\colon
\Fun(R\mathcal E)\longrightarrow\Fun(R\mathcal N)
\]
left adjoint to restriction from $R\mathcal N$ to
$R\mathcal E$.

The upper-plus construction defines a functor
\[
(-)^+\colon
\Fun(R\mathcal S)\longrightarrow\Fun(R\mathcal N)
\]
such that the composite with restriction to
$R\mathcal E$ is right adjoint to restriction from
$R\mathcal E$ to~$R\mathcal S$.

Both constructions preserve Green functors, and the
adjunctions restrict to the corresponding categories
of Green functors.
\end{Thm}
We also describe the induction and restriction maps and
the multiplicative structures explicitly in \cref{prop:lower-plus-action}, 
\cref{prop:upper-plus-action}, \cref{prop:lower-plus-product}, and 
\cref{con:upper-plus-product}.

As in the classical setting, the two adjunctions yield
a natural bijection
\[
\End^\gr_{R\mathcal S}
\bigl(\lRes{R\s}_{R\e}(A)
\bigr)
\cong
\Hom^\gr_{R\mathcal N}(A_+,A^+)
\]
for every Green $R\mathcal E$-functor $A$. The identity morphism
on the left corresponds to the mark morphism
$\mk^A\colon A_+\to A^+$.
\begin{Thm}(\Cref{sec:mark-map})
\label{thm:intro-mark-map}
For every Green $R\mathcal E$-functor $A$, the mark morphism
$\mk^A\colon A_+\to A^+$ is a morphism of Green
$R\mathcal N$-functors, natural in $A$.

If $|G|$ is invertible in $R$ for every
$G\in\mathfrak K$, then $\mk^A$ is an isomorphism.
Its inverse is given explicitly by M\"obius inversion.
\end{Thm}
Let $M$ be an $R\mathcal N$-functor. Using the Burnside ring action, 
we also establish an
induction--restriction decomposition of $M(G)$ when
$|G|$ is invertible in $R$. The decomposition yields a
zero-detection criterion used in the proof of canonical
induction; see
\Cref{prop:Burnside-decomposition,cor:Burnside-zero-detection}.

We now turn to canonical induction. Let $M$ be an
$R\mathcal N$-functor and $A$ an $R\mathcal E$-functor.
By lower-plus adjunction, an $R\mathcal E$-morphism
$\nu\colon A\to\iRes{R\e}{R\n}(M)$
determines a linearization morphism
$\lin^\nu\colon A_+\to M$ (\cref{con:linearization-map}).
By upper-plus adjunction, an $R\mathcal S$-morphism $\pi\colon\iRes{R\s}{R\n}(M)\longrightarrow\iRes{R\s}{R\e}(A)$
determines a tom Dieck morphism $\die^\pi\colon M\to A^+$ (\cref{con:tom-Dieck-map}).
When the mark morphism is invertible, the associated
canonical induction map (\cref{con:canonical-induction-map}) is
\[
\can^\pi=(\mk^A)^{-1}\circ\die^\pi
\colon M\longrightarrow A_+.
\]
Both $\die^\pi$ and $\can^\pi$ are
$R\mathcal E$-morphisms. We give an explicit formula
for $\can^\pi$ in \Cref{prop:canonical-induction-formula}.
When $M$ and $A$ are Green functors, we also
characterize its multiplicativity in terms of $\pi$ in
\Cref{prop:canonical-induction-multiplicative}.

We obtain the following analogue of Boltje's
canonical induction criterion.
\begin{Thm}(\Cref{sec:canonical-induction})
\label{thm:intro-canonical-induction}
Assume that $|G|$ is invertible in $R$ for every
$G\in\mathfrak K$. Let $M$ be an
$R\mathcal N$-functor and let $A$ be an
$R\mathcal E$-functor. Suppose that
$\nu\colon A\to\iRes{R\e}{R\n}(M)$
is an $R\mathcal E$-morphism and that
$\pi\colon
\iRes{R\s}{R\n}(M)
\longrightarrow
\iRes{R\s}{R\e}(A)$
is an $R\mathcal S$-morphism.
Then the following conditions are equivalent:
\begin{enumerate}[\rm(a)]
\item
$\lin^\nu\circ\can^\pi
=\id_{\iRes{R\e}{R\n}(M)}$.
\item
For every $G\in\mathfrak K$ and $m\in M(G)$,
\[
\nu_G(\pi_G(m))-m\in \sum\limits_{H<G}\lind{G}_H(M(H)).
\]
\end{enumerate}
\end{Thm}
%
In the classical setting of Mackey functors, the construction
and splitting criterion recover Boltje's results
\cite{Boltje98a}. Our formula specializes
over $\mathbb Q$ to the canonical induction formula
for $p$-permutation modules of \cite{BarkerMutlu19}.

Co\c{s}kun \cite{Coskun07} also applied the module-theoretic 
interpretation of plus constructions to the classification 
of simple Mackey functors. He subsequently used induction 
and coinduction to obtain alternative constructions of 
simple biset functors \cite{Coskun08}. In future work, 
we plan to investigate analogous questions for simple 
functors on Green categories and explore possible 
connections with the theory of simple modules over 
Green biset functors.

The paper is organized as follows.
\Cref{sec:preliminaries} contains the algebraic and
biset-theoretic preliminaries.
\Cref{sec:green-categories-from-mackey-systems}
constructs the Green category associated with a Mackey system,
and \Cref{sec:four-group-categories} develops the four-category
change-of-rings framework.
\Cref{sec:Burnside-idempotent-decomposition}
establishes the Burnside  action and associated idempotent decomposition.
\Cref{sec:plus-constructions} develops the plus constructions,
while \Cref{sec:mark-map} studies the mark morphism and its inverse.
Finally, \Cref{sec:canonical-induction} constructs the
canonical induction map and establishes the splitting criterion.

\begin{AI}
The author used generative AI tools as an aid in revising
this manuscript, including for language and organization.
All AI-assisted suggestions were critically reviewed by the author, who
takes full responsibility for the final manuscript, including all
mathematical statements and proofs.
\end{AI}

\begin{Ack}
This work grew out of research begun during the author's doctoral
studies at Bilkent University. The author is grateful to her doctoral
advisor, Laurence Barker, for his guidance during the development
of this work. The author would like to thank Paul Balmer for his 
valuable comments and Olcay Co\c{s}kun for bringing \cite{Coskun08} 
to her attention.
\end{Ack}

%
\section{Preliminaries}
\label{sec:preliminaries}
%
We collect the algebraic and biset-theoretic preliminaries used
throughout the paper. We first recall change-of-rings functors and
their adjunctions, then the required facts about bisets and generalized
biset morphisms, and finally the interpretation of group-category
functors as modules over locally unital quiver algebras.

%
\subsection{Induction, restriction and coinduction of modules}\
\label{ssec:induction-restriction-coinduction}
%

We first recall the change-of-rings functors and adjunctions
that will be used in the module-theoretic constructions below.
In this subsection, all rings are unital and all modules are
left modules unless otherwise stated. For a ring $\Gamma$,
we write $\Modcat(\Gamma)$ for the category of left
$\Gamma$-modules and $\Gamma$-linear maps.
\begin{Rec}
\label{Rec:change-of-rings}
Let $\alpha\colon\Lambda\to\Gamma$ be a unital ring homomorphism.
\begin{enumerate}
\item Any $\Gamma$-module can be regarded as a $\Lambda$-module via
$\alpha$. This induces a functor
\[
\lRes{\Lambda}_{\Gamma}^{\alpha}
\colon
\Modcat(\Gamma)\longrightarrow\Modcat(\Lambda),
\]
called restriction.

\item Regarding $\Gamma$ as a right $\Lambda$-module via $\alpha$,
induction is defined by
\[
\lInd{\Gamma}_{\Lambda}^{\alpha}
:=
\Gamma\otimes_{\Lambda}-
\colon
\Modcat(\Lambda)\longrightarrow\Modcat(\Gamma).
\]
\item Regarding $\Gamma$ as a left $\Lambda$-module via $\alpha$,
coinduction is defined by
\[
\lCo{\Gamma}_{\Lambda}^{\alpha}
:=
\hm_{\Lambda}(\Gamma,-)
\colon
\Modcat(\Lambda)\longrightarrow\Modcat(\Gamma).
\]
\end{enumerate}
\end{Rec}
We recall the standard change-of-rings adjunctions, see
\cite[Section 2.8]{Benson91}.
\begin{Rec}[Change-of-rings adjunctions]
\label{Rec:change-of-rings-adjunctions}
Let $\Lambda$ and $\Gamma$ be unital rings, and let
$\alpha\colon\Lambda\to\Gamma$ be a unital ring homomorphism.
Consider a $\Gamma$-module $M$ and a $\Lambda$-module $N$.
Then:
\begin{itemize}

\item[(a)] (Adjointness of induction and restriction.)
Naturally in $M$ and $N$, there is an isomorphism
\[
\hm_{\Gamma}
\bigl(
\lInd{\Gamma}_{\Lambda}^{\alpha}(N),M
\bigr)
\cong
\hm_{\Lambda}
\bigl(
N,\lRes{\Lambda}_{\Gamma}^{\alpha}(M)
\bigr).
\]
Maps
$\theta\colon
\lInd{\Gamma}_{\Lambda}^{\alpha}(N)\to M$
and
$\phi\colon
N\to\lRes{\Lambda}_{\Gamma}^{\alpha}(M)$
correspond if and only if~$\theta(\gamma\otimes n)=\gamma\phi(n)$
for all $\gamma\in\Gamma$ and $n\in N$.

\item[(b)] (Adjointness of restriction and coinduction.)
Naturally in $M$ and $N$, there is an isomorphism
\[
\hm_{\Lambda}
\bigl(
\lRes{\Lambda}_{\Gamma}^{\alpha}(M),N
\bigr)
\cong
\hm_{\Gamma}
\bigl(
M,\lCo{\Gamma}_{\Lambda}^{\alpha}(N)
\bigr).
\]
Maps
$\phi\colon
\lRes{\Lambda}_{\Gamma}^{\alpha}(M)\to N$
and
$\psi\colon
M\to\lCo{\Gamma}_{\Lambda}^{\alpha}(N)$
correspond if and only if
$\phi(\gamma m)=\psi(m)(\gamma)$
for all $\gamma\in\Gamma$ and $m\in M$.

\end{itemize}
\end{Rec}
\begin{Not}
\label{Not:change-of-rings}
When $\alpha$ is an inclusion
$\Lambda\hookrightarrow\Gamma$, we omit $\alpha$ from the notation
and write
\[
\iInd{\Gamma}{\Lambda},
\qquad
\iRes{\Lambda}{\Gamma},
\qquad
\iCo{\Gamma}{\Lambda}.
\]
When $\alpha$ is the quotient map
$\Lambda\twoheadrightarrow\Gamma=\Lambda/\Delta$
for an ideal $\Delta$ of $\Lambda$, we write
\[
\iDef{\Gamma}{\Lambda}
\qquad\text{and}\qquad
\iInf{\Lambda}{\Gamma}
\]
for induction and restriction, respectively.
\end{Not}
%

%
\subsection{Bisets and generalized biset morphisms}\
\label{ssec:bisets-generalized-bisets}
%

We next recall the biset notation and the generalized induction and
restriction morphisms needed to describe the Green categories used
later.
\begin{Rec}
\label{rec:bisets}
Let $F$ and $G$ be finite groups. We write
$B(F,G)=B(F\times G)$ for the Burnside group of finite $(F,G)$-bisets,
and denote the isomorphism class of an $(F,G)$-biset $X$ by $[X]$.
Every transitive $(F,G)$-biset is isomorphic to
$(F\times G)/L$
for some subgroup $L\leq F\times G$, and two transitive $(F,G)$-bisets
are isomorphic if and only if their point stabilizers are conjugate in
$F\times G$.

The biset category $\cC$ is the $\mathbb Z$-linear category whose objects
are finite groups and whose morphism group from $G$ to $F$ is
$\cC(F,G)=B(F,G)$;
see \cite[Chapters 2, 3]{Bouc10}. Thus an $(F,G)$-biset $X$ represents
the morphism $[X]\colon G\to F$ in $\cC$. If $X$ is an $(F,G)$-biset
and $Y$ is a $(G,H)$-biset, then
$[X][Y]=[X\times_G Y]$,
where $X\times_G Y$ denotes the set of $G$-orbits of $X\times Y$
under the action
$(x,y)\cdot g=(xg,g^{-1}y)$.

For a commutative unital ring $R$, we write $R\cC$ for the
$R$-linear category obtained from $\cC$ by extension of scalars:
$R\cC(F,G)=R\otimes_{\mathbb Z}B(F,G)$,
with composition extended $R$-bilinearly.
\end{Rec}
\begin{Def}[Generalized induction and restriction]
	\label{def:generalized-induction-restriction}
Given a group homomorphism $\alpha : G \rightarrow F,$ we define
transitive morphisms 
\[\lind{F}_{G}^{\alpha}=
\Bigg[ \frac{F\times G}{\{(\alpha(g),g): g \in G\}}\Bigg]
,\qquad\lres{G}_{F}^{\alpha}= \Bigg[ \frac{G\times F}
{\{(g, \alpha(g)) : g \in G\}}\Bigg]
\]
called induction and restriction, respectively. 
\end{Def}
\begin{Not}
\label{not:basic-biset-morphisms}
Following \cite{Barker16a}, we state basic results and set some notation. Let $\alpha :G \rightarrow F$ and
$\beta:H \rightarrow G $ be group homomorphisms. Then
\[
\lind{F}_G^\alpha\!\lind{}_H^\beta
=
\lind{F}_H^{\alpha\beta}
\qquad\text{and}\qquad
\lres{H}_G^\beta\!\lres{}_F^\alpha
=
\lres{H}_F^{\alpha\beta}.
\]
If $\alpha : G \hookrightarrow F$ is injective, $\lind{F}_G^\alpha$ is called
ordinary induction and $\lres{G}_F^\alpha$ is called ordinary
restriction. In this case 
we omit the symbol $\alpha$ from the notation, just
writing $\lind{F}_G$ and $\lres{G}_F$. If $\alpha$ is surjective, we have
\[
\ldef{F}_G
=
\lind{F}_G^\alpha
\qquad\text{and}\qquad
\linf{G}_F
=
\lres{G}_F^\alpha.
\]
For arbitrary $\alpha$, we have the factorizations
\[
\lind{F}_G^\alpha
=
\lind{F}_{\alpha(G)}\!
\deff_G
\qquad\text{and}\qquad
\lres{G}_F^\alpha
=
\linf{G}_{\alpha(G)}\!
\lres{}_F.
\]
If $\alpha$ is an isomorphism, we set
$\liso{F}_G^\alpha
=
\lind{F}_G^\alpha
=
\lres{F}_G^{\alpha^{-1}}$. 
If $g\in G$, then the isogation associated with the inner
automorphism $c_g$ is the identity morphism of $G$:
$\liso{G}_G^{c_g}=\id_G$.
\end{Not}
\begin{Rec}[Generalized Mackey relation]
\label{rec:generalized-mackey-relation}
We recall the generalized Mackey relation
\cite[Prop.~2.3]{Barker16b}. For a given diagram 
$F \stackrel{\alpha}\leftarrow I 
\stackrel{\theta}\rightarrow G \stackrel{\beta}
\leftarrow J \stackrel{\phi}\rightarrow H$
in the category of finite groups, we have
\[
\lind{F}_{I}^{\alpha}\!\res_{G}^{\theta}
\ind_{J}^{\beta}\!\res_{H}^{\phi}
=
\sum_{g\in\theta(I)\backslash G/\beta(J)}
\lind{F}_{L_g}^{\alpha\gamma_g}\!
\res_H^{\phi\psi_g}.
\]
where $I\stackrel{\gamma_{g}}\leftarrow L_{g} \stackrel{\psi_{g}}\rightarrow J$ 
is any pullback of 
$I \stackrel{\theta}\rightarrow G \stackrel{^{g}\beta}\leftarrow J$ where ${}^g\beta=c_g\circ\beta$.
\end{Rec}
\begin{Rec}[Right-free transitive morphisms]
\label{rec:right-free-transitive-morphisms}
An $(F,G)$-biset is called \emph{right-free} if its right
$G$-action is free. The right-free transitive morphisms
$F\leftarrow G$ are precisely the morphisms that can be expressed
in the form
\[
\lind{F}_V\!\lres{}_G^\beta
=
\lind{F}_V\!\inf_{\beta(V)}^\beta\res_G
=
\left[
\frac{F\times G}{\mathfrak G(V,\beta)}
\right].
\]
where $V\leq F$, $\beta\colon V\to G$, and
$\mathfrak G(V,\beta)
=
\{(v,\beta(v))\mid v\in V\}
\leq F\times G$.
\end{Rec}
\begin{Lem}[Uniqueness of right-free transitive presentations]
\label{lem:right-free-transitive-uniqueness}
Let $A,A'\leq G$, let $V,V'\leq F$, and let
$\rho:A\twoheadrightarrow V$ and $\rho':A'\twoheadrightarrow V'$.
Suppose that
\[
\lind{G}_A\res^{\rho}_V\res_F
=
\lind{G}_{A'}\res^{\rho'}_{V'}\res_F .
\]
Then there exist $g\in G$ and $f\in F$ such that
\[
A'={}^{g}A,\qquad
V'={}^{f}V,\qquad
\rho'\circ c_g=c_f\circ\rho .
\]
\end{Lem}
\begin{proof}
The two transitive bisets have stabilizers
$\mathfrak G(A,\rho)$ and $\mathfrak G(A',\rho')$ in $G\times F$.
By \Cref{rec:bisets}, these stabilizers are conjugate. Hence, for some
$g\in G$ and $f\in F$,
\[
\mathfrak G(A',\rho')
={}^{(g,f)}\mathfrak G(A,\rho).
\]
It follows that $A'={}^{g}A$ and
$\rho'c_g=c_f\rho$. Since $\rho$ and $\rho'$ are surjective,
$V'={}^{f}V$.
\end{proof}
%

%
\subsection{Group categories, functors and modules}\
\label{ssec:group-categories}
%

Throughout this subsection, let $R$ be a commutative unital ring.
Let $\mathcal L$ be a small $R$-linear category whose objects are
finite groups. By an $\mathcal L$-functor we mean an $R$-linear
functor $\mathcal L\longrightarrow R\text{-Mod}$.
The $\mathcal L$-functors can be regarded as modules over a locally
unital algebra via the following construction. We write $\Fun(\mathcal L)$ for the category of
$\mathcal L$-functors and $R$-linear natural transformations. When
$\Obj(\mathcal L)$ is infinite, the corresponding quiver algebra need
not be unital; its finite sums of identity morphisms instead form a
family of local units.
\begin{Rec}
\label{rec:locally-unital-algebras}
An $R$-algebra $\Lambda$, not necessarily unital, is called
\emph{locally unital} if every finite subset of $\Lambda$ has a
local unit, that is, an idempotent $u\in\Lambda$ such that
$u\lambda=\lambda=\lambda u$ for every element $\lambda$ of $\Lambda$.

A left $\Lambda$-module $M$ is called \emph{locally unital} if every
finite subset of $M$ is fixed by some idempotent of $\Lambda$.
For an arbitrary left $\Lambda$-module $M$, its \emph{locally unital
	part} is
\[
M_{\mathrm{lu}}
=
\{m\in M\mid um=m
\text{ for some idempotent }u\in\Lambda\}.
\]
Throughout the paper, modules over locally unital algebras are
understood to be locally unital.
\end{Rec}
\begin{Cons}[The quiver algebra]
\label{con:quiver-algebra}
For a small $R$-linear category $\mathcal L$, define its
\emph{quiver algebra} by
\[
{}^{\oplus}\mathcal L
=
\bigoplus_{F,G\in\Obj(\mathcal L)}
\mathcal L(F,G),
\]
where multiplication is given by composition when morphisms are
composable and is zero otherwise.

For every finite subset $\Omega\subseteq\Obj(\mathcal L)$, set
$u_\Omega=\sum_{G\in\Omega}\id_G$. These elements form a family of
local units for ${}^{\oplus}\mathcal L$, and every locally unital
${}^{\oplus}\mathcal L$-module $M$ decomposes as
\[
M=\bigoplus_{G\in\Obj(\mathcal L)}\id_GM.
\]
Thus, $\mathcal L$-functors are equivalently locally unital
${}^{\oplus}\mathcal L$-modules. Under this identification, an
$\mathcal L$-functor $M$ corresponds to
$\bigoplus_{G\in\Obj(\mathcal L)}M(G)$, where a morphism
$x\in\mathcal L(F,G)$ acts on the $G$-summand by
$M(x)\colon M(G)\to M(F)$ and annihilates the other summands.
Conversely, a locally unital ${}^{\oplus}\mathcal L$-module $M$
corresponds to the functor $G\mapsto\id_GM$.

We use this identification throughout when applying change-of-rings
constructions. In particular, for a group category $\mathcal X$, we
write $\Fun(R\mathcal X)$ for the category of $R\mathcal X$-functors,
identified with the category of locally unital modules over the
quiver algebra $R\mathcal X$.

For $R\mathcal X$-functors $A$ and $B$, we write
\[
\Hom_{R\mathcal X}(A,B)
:=
\Hom_{\Fun(R\mathcal X)}(A,B).
\]
\end{Cons}
\begin{Rec}[Locally unital change of rings]
\label{rec:change-of-rings-locally-unital}
Let $\Lambda$ and $\Gamma$ be locally unital $R$-algebras with a
common family of local units, and let
$\alpha\colon\Lambda\to\Gamma$
be an $R$-algebra homomorphism fixing every element of this family.
We consider only locally unital modules.

Restriction and induction are defined as in
\Cref{Rec:change-of-rings}. For coinduction, the usual module
$\hm_{\Lambda}(\Gamma,N)$ need not be locally unital, so we take its
locally unital part:
\[
\lCo{\Gamma}_{\Lambda}^{\alpha}(N)
=
\hm_{\Lambda}(\Gamma,N)_{\mathrm{lu}}.
\]
With these definitions, the change-of-rings adjunctions of
\Cref{Rec:change-of-rings-adjunctions} remain valid.
We use the notation of \Cref{Not:change-of-rings}.
\end{Rec}
\begin{Def}[Group categories]
\label{def:group-category}
Following \cite[Sec. 2]{Barker16b}, a \emph{group category}
$\mathcal X$ is a $\mathbb Z$-linear subcategory of the biset
category $\mathcal C$ such that, for all
$G,G'\in\Obj(\mathcal X)$, the transitive morphisms in
$\mathcal X(G',G)$ form a $\mathbb Z$-basis of
$\mathcal X(G',G)$.
\end{Def}
For a commutative unital ring $R$, we write $R\mathcal X$ for the
$R$-linear category obtained from $\mathcal X$ by extension of scalars,
so that
$R\mathcal X(G',G)
=
R\otimes_{\mathbb Z}\mathcal X(G',G)$.
By abuse of notation, we also write $R\mathcal X$ for its quiver
algebra ${}^{\oplus}(R\mathcal X)$.
\begin{Not}
\label{not:transitive-morphisms}
For a group category $\mathcal X$, let $\tm_{\mathcal X}$ denote
the set of its transitive morphisms, and set
$\tm_{\mathcal X}(G',G)
:=
\mathcal X(G',G)\cap\tm_{\mathcal X}$.
Thus, $\tm_{\mathcal X}(G',G)$ is an $R$-basis of
$R\mathcal X(G',G)$, and $\tm_{\mathcal X}$ is an $R$-basis of
$R\mathcal X$. In particular,
\[
R\mathcal X
=
\bigoplus_{\alpha\in\tm_{\mathcal X}}R\alpha.
\]
\end{Not}
%

%
\section{Green categories from Mackey systems}
\label{sec:green-categories-from-mackey-systems}
%

In this section, we recall the notions of Mackey systems and Green
categories that underlie our construction. Starting from a Mackey
system $\mathcal R$, we show that its injective morphisms form a
Mackey system $\mathcal I$ and that the pair
$(\mathcal I,\mathcal R)$ satisfies the Green-category axioms.
This framework will be applied in \Cref{con:four-group-categories}
to construct the ambient category $\mathcal N$.

Let $\mathfrak K$ be a set of finite groups closed under taking subgroups.
\begin{Rec}[Mackey systems]
	\label{rec:mackey-systems}
Following \cite[Sec.~3]{Barker16a}, a \emph{Mackey system} on
$\mathfrak K$ is a subcategory $\mathcal F$ of the category of finite
groups such that 
$\Obj(\mathcal{F}) = \mathfrak{K}$ and, writing
$\mathcal{F}(F,G)$ for the set of morphisms to a group 
$F \in \mathfrak{K}$ from a group $G \in \mathfrak{K}$,
the following four axioms hold:

\begin{enumerate}[\rm({MS}\,1)]
\item 
\label{it:MS-1}%
For every $I \leq G \in \mathfrak{K}$, 
the inclusion $I \hookrightarrow G$ is in $\mathcal{F}(G, I)$.
\item 
\label{it:MS-2}%
For every $I \leq G \in \mathfrak{K}$ and $g\in G$, 
the conjugation map $i\mapsto {}^{g}i$ belongs to $\mathcal{F}(^{g}I, I)$.
\item 
\label{it:MS-3}%
For every $\theta\in\mathcal F(F,G)$, the
corestriction $\ol{\theta}\colon G \to \theta(G)$
belongs to $\mathcal F(\theta(G),G)$.
\item 
\label{it:MS-4}%
If $\theta\in\mathcal F(F,G)$ is a group
isomorphism, then $\theta^{-1}\in\mathcal F(G,F)$.
\end{enumerate}
\end{Rec}
We shall use the injective morphisms in a Mackey system as a 
Mackey system in their own right, so we record the following elementary observation.
\begin{Lem}
\label{lem:injective-morphisms-mackey-system}
Let $\mathcal R$ be a Mackey system on $\mathfrak K$, and let
$\mathcal I$ be the subcategory of $\mathcal R$ consisting of its
injective morphisms. Then $\mathcal I$ is a Mackey system.
\end{Lem}
\begin{proof} 
This follows immediately from the axioms for a Mackey system. 
\end{proof}
\begin{Rem}
\label{rem:corestrictions}
Let $\theta\colon H\to G$ belong to $\Mor(\mathcal R)$, and let
$\ol{\theta}\colon H\onto\theta(H)$ be its corestriction.
By \MS{3}, $\ol{\theta}\in\Mor(\mathcal R)$.
If $\theta$ is injective, then $\ol{\theta}$ is a group isomorphism,
and hence an isomorphism in $\mathcal R$ by \MS{4}.
\end{Rem}
We next recall Barker's construction of Green categories from pairs
of Mackey systems.
\begin{Rec}[Green categories]
\label{rec:green-categories}
Following \cite[Sec.~2]{Barker16b}, consider a pair of Mackey systems
$(\mathcal I,\mathcal R)$ on $\mathfrak K$ such that
\begin{enumerate}[\rm({GC}\,1)]
\item 
\label{it:GC-1}%
The category $\mathcal{I}$ is a subcategory of $\mathcal{R}.$
\item 
\label{it:GC-2}%
Let $F \stackrel{\alpha}\leftarrow I \stackrel{\theta}\rightarrow G $
be a diagram with $\alpha \in \Mor(\mathcal{I})$ and $\theta \in\Mor(\mathcal{R})$.
If $\theta$ is surjective and there exists a group homomorphism $\beta$ such 
that $\alpha= \beta\theta$, then $\beta \in\Mor(\mathcal{I})$. If $\alpha$ 
is surjective and there exists a group homomorphism $\phi$ such that 
$\theta= \phi\alpha,$ then $\phi \in\Mor(\mathcal{R}).$
\item 
\label{it:GC-3}%
For every diagram
$I \stackrel{\theta} \rightarrow  G \stackrel{\beta} \leftarrow J $
with $\theta \in \Mor(\mathcal{R})$
and $\beta \in \Mor(\mathcal{I})$, there is a pullback  
$I \stackrel{\gamma}\leftarrow L \stackrel{\psi}\rightarrow J $ with $\gamma \in
\Mor(\mathcal{I})$ and $\psi \in\Mor(\mathcal{R})$.
\end{enumerate}

Let $(\mathcal I,\mathcal R)$ be a pair of Mackey systems satisfying 
the Green-category axioms. Using the induction and restriction morphisms 
of Definition~\ref{def:generalized-induction-restriction}, let $\mathcal G$ 
be the $\mathbb Z$-linear subcategory of $\mathcal C$ with 
$\Obj(\mathcal G)=\mathfrak K$ whose morphism groups are generated, as 
abelian groups, by finite compositions of morphisms $\ind^\alpha$ and $\res^\theta$, 
where $\alpha\in\Mor(\mathcal I)$ and $\theta\in\Mor(\mathcal R)$.

We call $\mathcal G$ the \emph{Green category} on $\mathfrak K$ associated with 
$(\mathcal I,\mathcal R)$; see \cite[Sec.~2]{Barker16b}. The morphisms in 
$\mathcal I$ are called the \emph{induction homomorphisms} for $\mathcal G$, and 
those in $\mathcal R$ are called the \emph{restriction homomorphisms} for $\mathcal G$.
\end{Rec}
We recall the following result from
\cite[Prop.~2.4]{Barker16b}.
\begin{Rec}[Transitive morphisms in a Green category]
\label{rec:green-category-transitive-morphisms}
Let $\mathcal G$ be the Green category associated with
$(\mathcal I,\mathcal R)$. Then $\mathcal G$ is a group category
in the sense of Definition~\ref{def:group-category}, and its
transitive morphisms are precisely the morphisms
$\ind^\alpha\res^\theta$,
where $\alpha\in\Mor(\mathcal I)$ and
$\theta\in\Mor(\mathcal R)$.
\end{Rec}
Barker~\cite[Ex.~2.6]{Barker16b} considers the case in which
$\mathcal R$ consists of all group homomorphisms and $\mathcal I$
of the injective ones. We now show that the same construction works
for an arbitrary Mackey system $\mathcal R$ by taking $\mathcal I$
to be its subcategory of injective morphisms.
\begin{Thm}[Green category from a Mackey system]
\label{thm:green-category-from-mackey-system}
Let $\mathcal R$ be a Mackey system on $\mathfrak K$, and let
$\mathcal I$ be its subcategory of injective morphisms.
Then $(\mathcal I,\mathcal R)$ satisfies \GC{1}--\GC{3}.
Consequently, the associated category of
\Cref{rec:green-categories} is a Green category.
\end{Thm}
\begin{proof}
By \Cref{lem:injective-morphisms-mackey-system}, $\mathcal I$ is a Mackey system, and \GC{1} is immediate.

For \GC{2}, let
$F \xleftarrow{\alpha} I \xrightarrow{\theta} G$
with $\alpha\in\Mor(\mathcal I)$ and $\theta\in\Mor(\mathcal R)$.
If $\theta$ is surjective and $\alpha=\beta\theta$, then $\theta$ is injective since $\alpha$ is injective, hence $\theta$ is an isomorphism. Thus
$\beta=\alpha\theta^{-1}\in\Mor(\mathcal R)$ by \MS{4}, and $\beta$ is injective, so $\beta\in\Mor(\mathcal I)$.
The other assertion is similar: if $\alpha$ is surjective, then $\alpha$ is an isomorphism and
$\varphi=\theta\alpha^{-1}\in\Mor(\mathcal R)$.

For \GC{3}, given
$I\xrightarrow{\theta}G\xleftarrow{\beta}J$
with $\theta\in\Mor(\mathcal R)$ and $\beta\in\Mor(\mathcal I)$, let
$L=\theta^{-1}(\beta(J))$.
The inclusion $\gamma:L\hookrightarrow I$ lies in $\Mor(\mathcal I)$ by \MS{1}, and there is a unique
$\psi:L\to J$ with $\theta\gamma=\beta\psi$. Since $\theta\gamma\in\Mor(\mathcal R)$, \MS{3} gives
$\overline{\theta\gamma}\colon L\to\theta(L)$ in
$\Mor(\mathcal R)$. Likewise,
$\overline\beta\colon J\to\beta(J)$ belongs to
$\Mor(\mathcal R)$ and is an isomorphism, so
$\overline\beta^{-1}\in\Mor(\mathcal R)$ by \MS{4}.
Since $\theta(L)\leq\beta(J)$, \MS{1} gives the inclusion
$\theta(L)\hookrightarrow\beta(J)$. Hence
\[
\psi=
\overline\beta^{-1}\circ
\bigl(\theta(L)\hookrightarrow\beta(J)\bigr)\circ
\overline{\theta\gamma}
\in\Mor(\mathcal R).
\]
By the definition of $L$ and the injectivity of $\beta$, the resulting
square is a pullback, proving \GC{3}.

Therefore $(\mathcal I,\mathcal R)$ satisfies \GC{1}--\GC{3}, and the associated category $\mathcal G$ is a Green category.
\end{proof}
\begin{Prop}[The injective induction system]
\label{prop:injective-induction-forced}
Let $(\mathcal I,\mathcal R)$ be a pair of Mackey systems satisfying
\GC{1}--\GC{3}. Suppose that every morphism in $\mathcal I$ is
injective. Then $\mathcal I$ is precisely the subcategory of
$\mathcal R$ consisting of its injective morphisms.
\end{Prop}
\begin{proof}
By \GC{1}, every morphism in $\mathcal I$ belongs to $\mathcal R$,
and by hypothesis it is injective. Thus $\mathcal I$ is contained in
the subcategory of injective morphisms of $\mathcal R$.

Conversely, let $\theta\colon H\to G$ be an injective morphism in
$\mathcal R$, and let
$\ol{\theta}\colon H\isoto\theta(H)$ be its corestriction. By
\MS{3}, $\ol{\theta}$ belongs to $\mathcal R$, and, since $\theta$
is injective, $\ol{\theta}$ is an isomorphism.

Apply \GC{2} to the diagram
$H\xleftarrow{\id_H}H\xrightarrow{\ol{\theta}}\theta(H)$.
Since $\ol{\theta}$ is surjective and
$\id_H=\ol{\theta}^{-1}\ol{\theta}$, it follows that
$\ol{\theta}^{-1}\in\Mor(\mathcal I)$. Hence
$\ol{\theta}\in\Mor(\mathcal I)$ by \MS{4} for $\mathcal I$.
The inclusion $\theta(H)\hookrightarrow G$ belongs to
$\mathcal I$ by \MS{1}. Thus
$\theta$, being the composite
$H\xrightarrow{\ol{\theta}}\theta(H)\hookrightarrow G$,
belongs to $\mathcal I$.
\end{proof}
%

%
\section{A quadruple of group categories}
\label{sec:four-group-categories}
%
The purpose of this section is to introduce the four group
categories $\n$, $\w$, $\e$, and $\s$ used in the module-theoretic description of the
plus constructions.  The main result of the section is the bimodule decomposition
$\mathcal W\otimes_{\mathcal S}\mathcal E\cong\mathcal N$.
This decomposition provides the change-of-rings framework
for the plus constructions and their adjunctions
established in \Cref{sec:plus-constructions}.

Let $\mathcal R$ be a Mackey system on $\mathfrak K$, and let
$\mathcal I$ be its subcategory of injective morphisms, as in
\Cref{lem:injective-morphisms-mackey-system}.

%
\subsection{The group categories $\mathcal S,\mathcal W,\mathcal E,\mathcal N$}\
\label{ssec:four-group-categories}
%

Following Co\c{s}kun's \cite{Coskun07} approach, we consider four group categories
$\mathcal S,\mathcal W,\mathcal E,\mathcal N$ on $\mathfrak K$.
\begin{Cons}[The four group categories]
\label{con:four-group-categories}
We define four $\mathbb Z$-linear subcategories
$\s$, $\w$, $\e$, $\n$ of the biset category
$\mathcal C$, all with object set $\mathfrak K$, as follows.

\begin{itemize}
\item $\mathcal N$ is generated by the morphisms
$\lind{G}_H^\alpha\res_F^\theta$,
where $\alpha\colon H\to G$ belongs to $\Mor(\mathcal I)$
and $\theta\colon H\to F$ belongs to $\Mor(\mathcal R)$.
	
\item $\mathcal W$ is generated by the morphisms
$\lind{G}_H^\alpha\res_F^\sigma$, 
where $\alpha\colon H\to G$ belongs to $\Mor(\mathcal I)$
and $\sigma\colon H\to F$ is a surjective morphism in
$\Mor(\mathcal R)$.
	
\item $\mathcal E$ is generated by the morphisms
$\lres{H}_F^\theta$, where $\theta\colon H\to F$ belongs to $\Mor(\mathcal R)$.
	
\item $\mathcal S$ is generated by the morphisms
$\lres{H}_F^\sigma$, where $\sigma\colon H\to F$ is a surjective morphism in
$\Mor(\mathcal R)$.
\end{itemize}

By Notation~\ref{not:basic-biset-morphisms}, the morphisms
generating $\mathcal S$ are inflation morphisms associated with
the surjective morphisms of $\mathcal R$.
\end{Cons}
Thus $\mathcal E$ contains all restriction morphisms, while
$\mathcal S$ contains only those arising from surjective
homomorphisms; similarly, $\mathcal W$ is obtained from
$\mathcal S$ by adjoining induction morphisms.
The resulting Green category $\mathcal N$ agrees with the inflation
Mackey category associated with $\mathcal R$ in the terminology of
\cite[Sec.~3]{Barker16a}.
\begin{Lem}
\label{lem:four-group-categories}
The categories
$\mathcal S,\mathcal W,\mathcal E,\mathcal N$
of \Cref{con:four-group-categories}
are group categories. Moreover,
$\mathcal S\subseteq\mathcal W\subseteq\mathcal N$,
$\mathcal S\subseteq\mathcal E\subseteq\mathcal N$,
and $\mathcal N$ is the Green category associated with
$(\mathcal I,\mathcal R)$.
\end{Lem}
\begin{proof}
By \Cref{con:four-group-categories,thm:green-category-from-mackey-system},
$\mathcal N$ is the Green category associated with
$(\mathcal I,\mathcal R)$.
Hence, by \Cref{rec:green-category-transitive-morphisms},
it is a group category, with transitive morphisms
$\lind{}^\alpha\res^\theta$ for
$\alpha\in\Mor(\mathcal I)$ and
$\theta\in\Mor(\mathcal R)$.

For $\mathcal E$ and $\mathcal S$, the specified transitive morphisms
are closed under composition by
\Cref{not:basic-biset-morphisms}; for $\mathcal S$, we also use that
composites of surjective homomorphisms are surjective.

For $\mathcal W$, let
$\lind{F}_I^\alpha\res_G^\sigma$ and
$\lind{G}_J^\beta\res_H^\tau$
be composable transitive generators.
By \Cref{rec:generalized-mackey-relation} and \GC{3}, their composite
is a sum of morphisms
$\lind{F}_{L_g}^{\alpha\gamma_g}\res_H^{\tau\psi_g}$,
where $\gamma_g\in\Mor(\mathcal I)$ and
$\psi_g\in\Mor(\mathcal R)$.
Since $\mathcal I$ and $\mathcal R$ are categories,
$\alpha\gamma_g\in\Mor(\mathcal I)$ and
$\tau\psi_g\in\Mor(\mathcal R)$. Moreover, $\psi_g$ is surjective,
being the pullback of the surjection $\sigma$, and hence
$\tau\psi_g$ is surjective.
Thus every summand is again a transitive morphism of $\mathcal W$.

Hence, in each case, the specified transitive morphisms span the
morphism groups. Since transitive bisets form a $\mathbb Z$-basis of
the morphism groups of $\mathcal C$, the specified transitive
morphisms form a basis. Thus
$\mathcal S,\mathcal W,\mathcal E,\mathcal N$ are group categories.

Finally,
$\lres{H}_F^\sigma
=
\lind{H}_H^{\id_H}\!\lres{}_F^\sigma$
gives $\mathcal S\subseteq\mathcal W$, and the remaining inclusions
are immediate from \Cref{con:four-group-categories}.
\end{proof}
\begin{Rem}[Two limiting cases]
\label{rem:four-group-categories-examples}
The construction of \Cref{con:four-group-categories} includes two
useful limiting cases.
If every morphism of $\mathcal R$ is injective, then
$\mathcal I=\mathcal R$, and the ambient category $\mathcal N$ is the
corresponding ordinary Mackey-type Green category.
At the other extreme, if $\mathcal R$ consists of all homomorphisms
between the groups in $\mathfrak K$, then $\mathcal I$ consists of all
injective homomorphisms, and $\mathcal N$ is Barker's category of
right-free bisets; see \cite[Ex.~2.6]{Barker16b}.

Thus, varying the Mackey system $\mathcal R$ varies the generalized
restriction morphisms available in $\mathcal E$, and hence the ambient
right-free Green category $\mathcal N$.
\end{Rem}

%
\subsection{Factorizations and the $\mathcal W\otimes_{\mathcal S}\mathcal E$ decomposition}\
\label{ssec:WE-decomposition}
%

The four categories are related by factorization of their transitive
morphisms. We first record two such factorizations and then use the
second one to identify $\mathcal N$ with
$\mathcal W\otimes_{\mathcal S}\mathcal E$.
\begin{Lem}[$\mathcal N$--$\mathcal E$ factorization]
\label{lem:N-E-factorization}
Let $F,G\in\mathfrak K$ and let
$\eta\in\tm_{\mathcal N}(G,F)$.
Then there exist a subgroup $U\leq G$ and a transitive morphism
$e\in\tm_{\mathcal E}(U,F)$ such that
$\eta=\lind{G}_Ue$.
Moreover, if
$\eta
=
\lind{G}_Ue
=
\lind{G}_Vf$
with
$e\in\tm_{\mathcal E}(U,F)$ and
$f\in\tm_{\mathcal E}(V,F)$,
then there exists $g\in G$ such that
$V={}^gU$ and $f=
\liso{V}_U^{c_g}e$.
\end{Lem}
\begin{proof}
By \Cref{rec:green-category-transitive-morphisms}, we may write
$\eta=\lind{G}_H^\alpha\res_F^\theta$
with
$\alpha\in\Mor(\mathcal I)$ and
$\theta\in\Mor(\mathcal R)$.
Since $\alpha$ is injective, its corestriction
$\ol{\alpha}\colon H\isoto\alpha(H)$
is an isomorphism in $\mathcal R$ by \MS{4}.
Hence, by \Cref{not:basic-biset-morphisms},
$\eta
=
\lind{G}_{\alpha(H)}
\res_F^{\theta\ol{\alpha}^{-1}}$. 
Since $\theta\ol{\alpha}^{-1}\in\Mor(\mathcal R)$, the second factor
belongs to $\tm_{\mathcal E}(\alpha(H),F)$ by
\Cref{con:four-group-categories}. This proves the first assertion.

Now suppose that
$\eta=\lind{G}_Ue=\lind{G}_Vf$.
Write $e=\lres{U}_F^\beta$ and
$f=\lres{V}_F^\gamma$, with
$\beta,\gamma\in\Mor(\mathcal R)$.
By \Cref{rec:right-free-transitive-morphisms,rec:bisets},
there exist $g\in G$ and $y\in F$ such that
$V={}^gU$ and $\gamma c_g=c_y\beta$.
Since inner isogations are identity morphisms by
\Cref{not:basic-biset-morphisms}, it follows that
$f=\liso{V}_U^{c_g}e$.
\end{proof}
The following decomposition is the analogue of
Co\c{s}kun's Mackey structure theorem
\cite[Thm.~3.2]{Coskun07}. We first record the factorization
of transitive morphisms that will be used in its proof.
\begin{Not}
\label{not:WE-quiver-algebras}
For the remainder of this section, when forming tensor products
or considering algebra structures on group categories, we regard
the corresponding $\mathbb Z$-linear categories as their locally
unital quiver algebras, as in
Construction~\ref{con:quiver-algebra}.
\end{Not}
\begin{Lem}[$\w$ --$\e$ factorization]
\label{lem:N-WE-factorization}
Let $\eta\in\tm_{\mathcal N}(G,F)$. Then there exist a subgroup
$V\leq F$ and a transitive morphism
$w\in\tm_{\mathcal W}(G,V)$ such that
$\eta=w\lres{V}_F$.

Moreover, if
$\eta
=
w\lres{V}_F
=
w'\lres{V'}_F$,
where $V,V'\leq F$,
$w\in\tm_{\mathcal W}(G,V)$, and
$w'\in\tm_{\mathcal W}(G,V')$,
then there is an isogation
$s\in\tm_{\mathcal S}(V',V)$
such that
$w=w's$ and
$s\lres{V}_F=\lres{V'}_F$.
Consequently,
\[
w\otimes_{\mathcal S}\lres{V}_F
=
w'\otimes_{\mathcal S}\lres{V'}_F.
\]
\end{Lem}
\begin{proof}
Write
$\eta=\lind{G}_H^\alpha\res_F^\theta$
as in \Cref{rec:green-category-transitive-morphisms}. Factoring
$\theta$ through its image gives
\[
\eta
=
\bigl(
\lind{G}_H^\alpha
\res_{\theta(H)}^{\ol{\theta}}
\bigr)
\lres{\theta(H)}_F,
\]
where the first factor belongs to
$\tm_{\mathcal W}(G,\theta(H))$.

For uniqueness, write
$w=\lind{G}_A\res_V^\rho$ and
$w'={}\lind{G}_{A'}\res_{V'}^{\rho'}$,
with $\rho\colon A\twoheadrightarrow V$ and
$\rho'\colon A'\twoheadrightarrow V'$.
By \Cref{lem:right-free-transitive-uniqueness}, there exist
$g\in G$ and $f\in F$ such that
$A'={}^{g}A$, $V'={}^{f}V$, and
$\rho'c_g=c_f\rho$.
Thus, for $s=\liso{V'}_V^{c_f}$,
we have
$w=w's$ and
$s\lres{V}_F=\lres{V'}_F$.
The final equality follows from the balancing relation.
\end{proof}
\begin{Prop}[The bimodule decomposition]
\label{prop:WE-decomposition}
Composition induces an isomorphism of
$(\mathcal W,\mathcal E)$-bimodules
\begin{equation}
\label{eq:WE-decomposition}
\Phi\colon
\mathcal W\otimes_{\mathcal S}\mathcal E
\isoto
\mathcal N,
\qquad
w\otimes_{\mathcal S}e\longmapsto we.
\end{equation}
\end{Prop}
\begin{proof}
By \Cref{lem:four-group-categories}, composition is
$\mathcal S$-balanced and hence induces the bimodule morphism
$\Phi$ of \Cref{eq:WE-decomposition}.
By \Cref{lem:N-WE-factorization}, every transitive morphism of
$\mathcal N$ lies in its image, so $\Phi$ is surjective.

For $\eta\in\tm_{\mathcal N}(G,F)$, choose a factorization
$\eta=w\lres{V}_F$ as in \Cref{lem:N-WE-factorization} and set
\[
\Psi(\eta)
=
w\otimes_{\mathcal S}\lres{V}_F.
\]
The uniqueness assertion of \Cref{lem:N-WE-factorization} shows
that this is well defined. Since the transitive morphisms form a
$\mathbb Z$-basis of $\mathcal N$, it extends linearly to
$\Psi\colon\mathcal N\to
\mathcal W\otimes_{\mathcal S}\mathcal E$.
By construction, $\Phi\Psi=\id_{\mathcal N}$.

It remains to check $\Psi\Phi=\id$.
Since the transitive morphisms span $\mathcal W$ and
$\mathcal E$, it is enough to consider
composable transitive morphisms
$w\otimes_{\mathcal S}\lres{K}_F^\theta$.
By \Cref{rem:corestrictions},
$\lres{K}_{\theta(K)}^{\ol{\theta}}\in\mathcal S$.
By \Cref{not:basic-biset-morphisms},
$w\lres{K}_{\theta(K)}^{\ol{\theta}}$
is again a transitive morphism of $\mathcal W$.
Hence
\[
\Psi\Phi
\bigl(
w\otimes_{\mathcal S}\lres{K}_F^\theta
\bigr)
=
w\lres{K}_{\theta(K)}^{\ol{\theta}}
\otimes_{\mathcal S}\lres{\theta(K)}_F
=
w\otimes_{\mathcal S}\lres{K}_F^\theta
\]
by the balancing relation and
\Cref{not:basic-biset-morphisms}.
Thus $\Psi\Phi=\id$, and the result follows.
\end{proof}
\begin{Cor}
\label{cor:WE-decomposition-scalar-extension}
For every commutative unital ring $R$, multiplication induces an
isomorphism of $(R\mathcal W,R\mathcal E)$-bimodules
\[
R\mathcal W\otimes_{R\mathcal S}R\mathcal E
\isoto
R\mathcal N.
\]
\end{Cor}
\begin{proof}
This follows from
Proposition~\ref{prop:WE-decomposition}
by extension of scalars.
\end{proof}
\begin{Rem}[Relation with admissible biset categories]
\label{rem:admissible-biset-categories}
Let
$\eta
=
\lind{G}_H^\alpha\res_F^\theta
\in\tm_{\mathcal N}(G,F),
$
with $\alpha\in\Mor(\mathcal I)$ and
$\theta\in\Mor(\mathcal R)$. A point stabilizer of the corresponding
transitive $(G,F)$-biset is
$D_{\alpha,\theta}
=
\{(\alpha(h),\theta(h))\mid h\in H\}
\leq G\times F$.
Writing $p_1$ and $p_2$ for the two projections and
$k_2(D)
=
\{f\in F\mid (1,f)\in D\}$,
we have
$k_2(D_{\alpha,\theta})=1$,
$p_1(D_{\alpha,\theta})=\alpha(H)$,
$p_2(D_{\alpha,\theta})=\theta(H)$.
Thus every transitive morphism of $\mathcal N$ is right-free.
Moreover,
\[
\begin{aligned}
	\eta\in\tm_{\mathcal E}(G,F)
	&\Longleftrightarrow p_1(D_{\alpha,\theta})=G,\\
	\eta\in\tm_{\mathcal W}(G,F)
	&\Longleftrightarrow p_2(D_{\alpha,\theta})=F,\\
	\eta\in\tm_{\mathcal S}(G,F)
	&\Longleftrightarrow
	p_1(D_{\alpha,\theta})=G
	\text{ and }
	p_2(D_{\alpha,\theta})=F.
\end{aligned}
\]

This places the quadruple
$(\mathcal S,\mathcal W,\mathcal E,\mathcal N)$ in the right-free
part of the more general admissible-biset framework of
\cite{BoltjeRaggiCardenasValeroElizondo19}.
Their framework allows more general choices of biset operations,
including deflations. However, their upper-plus construction
requires right-free bisets. Moreover, the mark morphism need
not commute with deflations in general, whereas its naturality
is established in the right-free setting; see
\cite[Sec.~5 and Thm.~6.1]{BoltjeRaggiCardenasValeroElizondo19}.
Our framework thus concerns Green categories with injective induction 
homomorphisms, whose transitive morphisms are right-free.
\end{Rem}

%
\subsection{Module-theoretic consequences}\
\label{ssec:module-theoretic-consequences}
%

Using the identification of $R\mathcal X$-functors with locally
unital modules from \Cref{con:quiver-algebra}, we now pass from the
categorical decomposition to the corresponding functor categories.
The subcategory $\mathcal S\subseteq\mathcal W$ admits a
complementary span generated by the transitive morphisms with
nontrivial induction part. We show that this complement is an ideal,
producing a projection $\mathcal W\to\mathcal S$. Together with the
decomposition of the preceding subsection, this provides the
restriction, induction, inflation, and coinduction maps used in the 
plus constructions of \Cref{sec:plus-constructions}.
\begin{Not}
\label{not:Delta-W}
Using Notation~\ref{not:transitive-morphisms}, we set
\[
\Delta_{\mathcal W}
:=
\bigoplus_{\omega\in
	\tm_{\mathcal W}\setminus\tm_{\mathcal S}}
\mathbb Z\omega.
\]
Since $\mathcal S\subseteq\mathcal W$ by
Lemma~\ref{lem:four-group-categories}, we have a decomposition
of $\mathbb Z$-modules
$\mathcal W
=
\mathcal S\oplus\Delta_{\mathcal W}$.
\end{Not}
\begin{Lem}
\label{lem:Delta-W-transitive}
Let
$\omega={}_G\ind_H^\alpha\res_F^\sigma$
be a transitive morphism in $\mathcal W$, where
$\alpha\in\Mor(\mathcal I)$ and
$\sigma\in\Mor(\mathcal R)$ is surjective.
Then $\omega\in\mathcal S$ if and only if $\alpha(H)=G$.
Consequently, $\Delta_{\mathcal W}$ is spanned by the
transitive morphisms of $\mathcal W$ whose induction part is
proper.
\end{Lem}
\begin{proof}
By \Cref{rem:corestrictions,not:basic-biset-morphisms},
$\omega
=
\lind{G}_{\alpha(H)}
\res_F^{\sigma\ol{\alpha}^{-1}}$,
where $\sigma\ol{\alpha}^{-1}$ is surjective.
Thus $\alpha(H)=G$ implies $\omega\in\mathcal S$.

Conversely, if $\omega\in\mathcal S$, then its point stabilizer has
first projection $\alpha(H)$, while any transitive expression of
$\omega$ in $\mathcal S(G,F)$ has first projection $G$.
By \Cref{rec:bisets}, these stabilizers are conjugate, so
$\alpha(H)=G$.

The final assertion follows from \Cref{not:Delta-W}.
\end{proof}
\begin{Prop}
\label{prop:Delta-W-ideal}
The $\mathbb Z$-submodule
$\Delta_{\mathcal W}$ is a two-sided ideal of
$\mathcal W$.
\end{Prop}
\begin{proof}
It is enough to consider composable transitive morphisms
$x=\lind{F}_I^\alpha\res_G^\sigma$ and
$y=\lind{G}_J^\beta\res_H^\tau$ in $\mathcal W$.
As in the proof of \Cref{lem:four-group-categories},
\Cref{rec:generalized-mackey-relation} expresses $xy$ as a sum of
transitive morphisms
$\lind{F}_{L_g}^{\alpha\gamma_g}\res_H^{\tau\psi_g}$,
where
$I\stackrel{\gamma_g}\leftarrow L_g
\stackrel{\psi_g}\rightarrow J$
is a pullback of
$I\stackrel{\sigma}\rightarrow G
\stackrel{{}^g\beta}\leftarrow J$.

If $x\in\Delta_{\mathcal W}$, then
$\alpha(I)<F$ by \Cref{lem:Delta-W-transitive}, and hence
$(\alpha\gamma_g)(L_g)<F$ for every $g$.
If $y\in\Delta_{\mathcal W}$, then
${}^g\beta(J)<G$. Since $\sigma\colon I\twoheadrightarrow G$ is
surjective and ${}^g\beta(J)<G$, its inverse image
$\sigma^{-1}({}^g\beta(J))=\gamma_g(L_g)$ is a proper subgroup of
$I$. Therefore
$(\alpha\gamma_g)(L_g)<F$.

Thus every summand belongs to $\Delta_{\mathcal W}$ by
\Cref{lem:Delta-W-transitive}.
\end{proof}
\begin{Cons}
\label{con:W-to-S-projection}
By Notation~\ref{not:Delta-W} and
Proposition~\ref{prop:Delta-W-ideal}, the projection
\[
\pi_{\mathcal S}^{\mathcal W}\colon
\mathcal W\longrightarrow\mathcal S
\]
with kernel $\Delta_{\mathcal W}$ is a homomorphism of locally
unital $\mathbb Z$-algebras.  Since $\pi_{\mathcal S}^{\mathcal W}$ restricts to the identity
on $\mathcal S$, it fixes every identity morphism and hence
preserves the common local units. After extension of scalars, we obtain
a homomorphism
$R\mathcal W\to R\mathcal S$ with kernel
$R\Delta_{\mathcal W}$.
\end{Cons}

Thus, we have a diagram of locally unital algebra homomorphisms
\[
\xymatrix@C=3ex@R=3ex{
	&& \mathcal N &\\
	& \mathcal W\ar@{ ->>}[ld]_{\pi_{\mathcal S}^{\mathcal W}}\ar@{ >->}[ur]
	&&\mathcal E\ar@{ >->}[ul]\\
	\mathcal S\ar@{=}[rr]&&\mathcal S\ar@{ >->}[ul]\ar@{ >->}[ur]&
}
\]
where the upward arrows are inclusions.

By Construction~\ref{con:quiver-algebra}, the four quiver algebras
$R\mathcal S$, $R\mathcal W$, $R\mathcal E$, and $R\mathcal N$
have the same local units, namely the finite sums of identity
morphisms.
\begin{Not}
\label{not:W-restriction-inflation}
The inclusions among the four group categories give restriction
functors
\[
\begin{aligned}
\iRes{R\s}{R\w}
&\colon
\Fun(R\w)\longrightarrow\Fun(R\s),
&
\iRes{R\s}{R\e}
&\colon
\Fun(R\e)\longrightarrow\Fun(R\s),
\\
\iRes{R\w}{R\n}
&\colon
\Fun(R\n)\longrightarrow\Fun(R\w),
&
\iRes{R\e}{R\n}
&\colon
\Fun(R\n)\longrightarrow\Fun(R\e).
\end{aligned}
\]
The projection
$R\mathcal W\twoheadrightarrow R\mathcal S$
of \Cref{con:W-to-S-projection} gives the inflation functor
\[
\iInf{R\w}{R\s}
\colon
\Fun(R\s)\longrightarrow\Fun(R\w).
\]
Thus, for every $R\s$-functor $B$,
$\iInf{R\w}{R\s}(B)$ has the same underlying
$R\s$-functor as $B$, with $R\Delta_{\mathcal W}$ acting trivially
under the module identification of \Cref{con:quiver-algebra}.
\end{Not}
\begin{Lem}
\label{lem:Res-Inf-identity}
There is a natural isomorphism
\[
{}\iRes{R\mathcal S}{R\mathcal W}\!\Inf_{R\mathcal S}
\cong
\id_{\Fun(R\s)}.
\]
\end{Lem}
\begin{proof}
This is immediate from \Cref{not:W-restriction-inflation}.
\end{proof}
The following comparison result is the analogue of
Co\c{s}kun's \cite[Thm.~3.4]{Coskun07}.
\begin{Prop}[Induction and coinduction comparisons]
\label{prop:Ind-Coind-comparison}
The following natural isomorphisms hold.
\begin{enumerate}[\rm(a)]
\item
$\iRes{R\w}{R\n}\!\iInd{}{R\e}
\cong
\iInd{R\w}{R\s}\!\iRes{}{R\e}$.
\item $\iCo{R\e}{R\s}\!
\iRes{}{R\w}
\cong
{}\iRes{R\e}{R\n}\!
\iCo{}{R\w}$.
\end{enumerate}

For an $R\mathcal W$-module $M$, under the isomorphism in
\textup{(b)}, elements
\[
\epsilon\in
{}\iCo{R\mathcal E}{R\mathcal S}\!
\iRes{}{R\mathcal W}(M)
\quad\text{and}\quad
\nu\in
{}\iRes{R\mathcal E}{R\mathcal N}\!
\iCo{}{R\mathcal W}(M)
\]
correspond if and only if
$\nu(we)=w\epsilon(e)$
for all $w\in R\mathcal W$ and $e\in R\mathcal E$.
Equivalently, $\epsilon$ is the restriction of $\nu$ to
$R\mathcal E$.
\end{Prop}
\begin{proof}
Part \textup{(a)} follows from
\Cref{cor:WE-decomposition-scalar-extension} by tensoring over
$R\e$.

For \textup{(b)}, we obtain natural $R\e$-isomorphisms
\[
\begin{aligned}
	\iRes{R\e}{R\n}\!
	\iCo{}{R\w}(M)
	&=
	\Hom_{R\w}
	\bigl(
	{}_{R\w}R\n_{R\e},M
	\bigr)_{\mathrm{lu}}
	&&\text{by \Cref{rec:change-of-rings-locally-unital}}
	\\[.1em]
	&\cong
	\Hom_{R\w}
	\bigl(
	{}_{R\w}
	R\w\otimes_{R\s}R\e_{R\e},
	M
	\bigr)_{\mathrm{lu}}
	&&\text{by \Cref{cor:WE-decomposition-scalar-extension}}
	\\[.1em]
	&\cong
	\Hom_{R\s}
	\bigl(
	R\e,
	\iRes{R\s}{R\w}(M)
	\bigr)_{\mathrm{lu}}
	&&\text{by \Cref{Rec:change-of-rings-adjunctions}}
	\\[.1em]
	&=
	\iCo{R\e}{R\s}\!
	\iRes{}{R\w}(M)
	&&\text{by \Cref{rec:change-of-rings-locally-unital}.}
\end{aligned}
\]
Under these isomorphisms, $\epsilon$ corresponds to $\nu$ precisely
when
$\nu(we)=w\epsilon(e)$.
Choosing a local unit $u\in R\s$ with $ue=e$, we obtain
\[
\nu(e)=\nu(ue)=u\epsilon(e)=\epsilon(e),
\]
so $\epsilon$ is the restriction of $\nu$ to $R\e$.
Conversely, if $\epsilon=\nu|_{R\e}$, then
$\nu(we)=w\nu(e)=w\epsilon(e)$ by $R\w$-linearity.
\end{proof}

%
\subsection{Green functors}\
\label{ssec:green-functors}
%

Following Barker~\cite[Sec.~3]{Barker16b}, we recall
the notion of a Green functor for the Green category $\n$;
see also Romero~\cite[Sec.~2]{Romero12}.
We also introduce Green $R\s$- and $R\e$-functors,
which will be used to establish the multiplicative structures
of the upper- and lower-plus constructions in
\Cref{sec:plus-constructions}.

\begin{Def}[Green $R\e$- and $R\s$-functors]
\label{def:Green-E-functor}
\label{def:Green-S-functor}
A \emph{Green $R\e$-functor} is an $R\e$-functor $A$ such that:
\begin{enumerate}
\item $A(G)$ is a unital $R$-algebra for every
$G\in\mathfrak K$;

\item for every $\theta\colon I\to G$ in
$\Mor(\mathcal R)$, the restriction map
$\lres{I}_G^\theta\colon A(G)\longrightarrow A(I)$
is a unital $R$-algebra homomorphism.
\end{enumerate}
A \emph{Green $R\s$-functor} is defined similarly and with condition~(b)
required only for surjective $\theta\in\Mor(\mathcal R)$.
\end{Def}
\begin{Def}[Green $R\mathcal N$-functors]
\label{def:Green-N-functor}
A \emph{Green $R\mathcal N$-functor} is an
$R\mathcal N$-functor $A$ such that:
\begin{enumerate}
\item $A(G)$ is a unital $R$-algebra for every
$G\in\mathfrak K$;

\item for every $\theta\colon I\to G$ in
$\Mor(\mathcal R)$, the restriction map
$\lres{I}_G^\theta\colon A(G)\to A(I)$
is a unital $R$-algebra homomorphism;

\item for every $\alpha\colon I\to F$ in $\Mor(\mathcal I)$,
$a\in A(F)$, and $b\in A(I)$, the Frobenius relations
\[
{}\lind{F}_I^\alpha
\bigl({}\lres{I}_F^\alpha(a)\cdot b\bigr)
=
a\cdot{}\lind{F}_I^\alpha(b)
\]
and
\[
{}\lind{F}_I^\alpha
\bigl(b\cdot{}\lres{I}_F^\alpha(a)\bigr)
=
{}\lind{F}_I^\alpha(b)\cdot a
\]
hold.
\end{enumerate}
\end{Def}
\begin{Def}[Morphisms of Green functors]
\label{def:Green-functor-morphisms}
Let $\mathcal X$ be $\s$, $\e$ or $\n$.
A morphism $f\colon A\to B$ of Green $R\mathcal X$-functors is
a morphism of $R\mathcal X$-functors such that, for every
$G\in\mathfrak K$, the component
$f_G\colon A(G)\to B(G)$ is a unital $R$-algebra homomorphism.
\end{Def}
\begin{Not}
\label{not:Green-functor-categories}
For $\mathcal X\in\{\s,\e,\n\}$, we write
\[
\grfun{R\mathcal X}
\]
for the category of Green $R\mathcal X$-functors and the morphisms of
\Cref{def:Green-functor-morphisms}.

For Green $R\mathcal X$-functors $A$ and $B$, we write
\[
\Hom_{R\mathcal X}^\gr(A,B)
:=
\Hom_{\grfun{R\mathcal X}}(A,B).
\]
\end{Not}
By Lemma~\ref{lem:four-group-categories},
Definition~\ref{def:Green-N-functor} is precisely Barker's
definition of a Green functor for the Green category $\mathcal N$.
\begin{Rem}
\label{rem:Green-restrictions}
Restriction induces functors
\[
\iRes{R\s}{R\e}
\colon
\grfun{R\e}\longrightarrow\grfun{R\s}
\quad \text{and}\quad \iRes{R\e}{R\n}
\colon
\grfun{R\n}\longrightarrow\grfun{R\e}. 
\]
Consequently,
$\iRes{R\s}{R\n}$
also restricts to a functor
$\grfun{R\n}\to\grfun{R\s}$.
\end{Rem}
%

%
\section{The Burnside idempotent decomposition}
\label{sec:Burnside-idempotent-decomposition}
%

The purpose of this section is to separate, for an $R\n$-functor $M$,
the part of $M(G)$ detected by restriction to proper subgroups from the
part generated by induction from proper subgroups. We do this using
the action of the Burnside ring on $M(G)$ and the primitive
Burnside idempotent $e_G^G$. When $|G|$ is invertible in $R$, this
yields an induction--restriction decomposition of $M(G)$ that will be
used in the canonical-induction criterion of
\Cref{sec:canonical-induction}.

%
\subsection{The Burnside action}\
\label{ssec:Burnside-action}
%

We begin by recalling the Burnside-ring notation that will be used
throughout this section.
\begin{Rec}
\label{rec:Burnside-ring}
Let $G\in\mathfrak K$.  We write $B(G)$ for the Burnside ring of~$G$
and set $d_U^G:=[G/U]$ for $U\leq G$.
The elements $d_U^G$, with $U$ running through representatives of the
$G$-conjugacy classes of subgroups of~$G$, form a $\mathbb Z$-basis
of~$B(G)$.  For $U,V\leq G$, multiplication is given by
\[
d_U^G\cdot d_V^G
=
\sum_{g\in U\backslash G/V}
d_{U\cap{}^gV}^{G}.
\]
If $H\leq G$, restriction is given by
\[
{}\ires{H}{G}(d_U^G)
=
\sum_{g\in H\backslash G/U}
d_{H\cap{}^gU}^{H}.
\]
We write $RB(G):=R\otimes_{\mathbb Z}B(G)$.
\end{Rec}
The Burnside algebra $RB(G)$ embeds naturally in the
endomorphism algebra of~$G$.
\begin{Prop}[The Burnside embedding]
\label{prop:Burnside-embedding}
For every $G\in\mathfrak K$, there is a unital injective
$R$-algebra homomorphism
$\Delta_G\colon RB(G)\longrightarrow R\n(G,G)$
determined by 
\[
\Delta_G(d_U^G)={}\lind{G}_U\!\lres{}_G
\]
for $U\leq G$.
\end{Prop}
\begin{proof}
If $U$ and $V$ are $G$-conjugate, then
${}\lind{G}_U\!\lres{}_G={}\lind{G}_V\!\lres{}_G$
by \Cref{rec:bisets}, so the assignment is well defined.
Moreover, these morphisms belong to $R\n(G,G)$ by \MS{1}
and the definition of $\mathcal I$.

By \Cref{rec:generalized-mackey-relation},
\[
{}\lind{G}_U\!\lres{}_G\!\lind{}_V\!\lres{}_G
=
\sum_{g\in U\backslash G/V}
{}\lind{G}_{{U\cap{}^gV}}\!\lres{}_G,
\]
so comparison with \Cref{rec:Burnside-ring} shows that
$\Delta_G$ is multiplicative. Moreover,
$\Delta_G(d_G^G)=\id_G$, so $\Delta_G$ is unital.

Finally, ${}\lind{G}_U\!\lres{}_G$ has point stabilizer
$\{(u,u)\mid u\in U\}\leq G\times G$.
Hence non-conjugate subgroups give distinct transitive morphisms by
\Cref{rec:bisets}, and these are $R$-linearly independent by
\Cref{not:transitive-morphisms}.
Thus $\Delta_G$ is injective.
\end{proof}
The resulting action agrees with the standard Burnside action on a
Mackey functor; compare \cite[Proposition~1.5(i)]{Boltje98a}.
\begin{Cons}[The Burnside action]
\label{con:Burnside-action}
Let $M$ be an $R\n$-functor and let $G\in\mathfrak K$.
Via \Cref{prop:Burnside-embedding}, we regard $M(G)$ as an
$RB(G)$-module by $x\cdot m:=\Delta_G(x)(m)$ for
$x\in RB(G)$ and $m\in M(G)$.

In particular, for $U\leq G$, we have
$d_U^G\cdot m={}\lind{G}_U\bigl(\lres{U}_G(m)\bigr)$.
\end{Cons}
\begin{Lem}[Compatibility of the Burnside action]
\label{lem:Burnside-action-compatibility}
Let $M$ be an $R\n$-functor and let $H\leq G$.
\begin{enumerate}
\item[(a)] For $x\in RB(G)$ and $m\in M(G)$, we have
\[
{}\lres{H}_G(x\cdot m)
=
\bigl({}\lres{H}_G(x)\bigr)\cdot\bigl(\lres{H}_G(m)\bigr).
\]
\item[(b)] For $x\in RB(G)$ and $n\in M(H)$, we have
\[
x\cdot{}\lind{G}_H(n)
=
{}\lind{G}_H
\left(
\bigl({}\lres{H}_G(x)\bigr)\cdot n
\right).
\]
\end{enumerate}
\end{Lem}
\begin{proof}
By $R$-linearity, it is enough to take $x=d_U^G$.
Both assertions follow from the restriction formula in
\Cref{rec:Burnside-ring} and the generalized Mackey relation
\Cref{rec:generalized-mackey-relation}; for~(b), reindex the
Mackey sum by $g\mapsto g^{-1}$.
\end{proof}

%
\subsection{Primitive idempotents and the induction--restriction decomposition}\
\label{ssec:Burnside-idempotent-decomposition}
%

We recall the primitive idempotents of the rational Burnside algebra;
see \cite[Remark~2.5]{Boltje98a}.
\begin{Rec}[Primitive Burnside idempotents]
\label{rec:primitive-Burnside-idempotents}
Let $G\in\mathfrak K$.  For $H\leq G$, set
\[
e_H^G
:=
\frac{1}{|N_G(H)|}
\sum_{K\leq H}
|K|\,\mu(K,H)\,d_K^G
\in
\mathbb Q\otimes_{\mathbb Z}B(G).
\]
When $H$ runs through representatives of the $G$-conjugacy
classes of subgroups of~$G$, the elements $e_H^G$ form a complete
set of primitive orthogonal idempotents of the rational Burnside
algebra $\mathbb Q\otimes_{\mathbb Z}B(G)$.

If $|G|$ is invertible in $R$, then so is $|N_G(H)|$ for every
$H\leq G$, and the same formulas define orthogonal idempotents
in $RB(G)$.  In particular,
\[
e_G^G
=
\frac{1}{|G|}
\sum_{K\leq G}
|K|\,\mu(K,G)\,d_K^G,
\]
and therefore
\[
1-e_G^G
=
-\frac{1}{|G|}
\sum_{K<G}
|K|\,\mu(K,G)\,d_K^G.
\]
Moreover, ${}_H\res_G(e_G^G)=0$ for every strict subgroup
$H<G$; see \cite[proof of Proposition~6.2]{Boltje98a}.
\end{Rec}
\begin{Not}
\label{not:induction-restriction-submodules}
Let $M$ be an $R\n$-functor and let $G\in\mathfrak K$.  
Following \cite[Subsec.~1.6--1.7]{Boltje98a}, we set
\[
\cI(M)(G)
:=
\sum_{H<G}
{}\lind{G}_H\bigl(M(H)\bigr)
=
\sum_{H<G}
\Img
\bigl({}\lind{G}_H\colon M(H)\to M(G)\bigr)
\]
and
\[
\cK(M)(G)
:=
\bigcap_{H<G}
\Ker
\bigl({}\lres{H}_G\colon M(G)\to M(H)\bigr).
\]
\end{Not}
The following decomposition is the analogue, for $R\n$-functors, of
\cite[Proposition~6.2]{Boltje98a}.
\begin{Prop}[The induction--restriction decomposition]
\label{prop:Burnside-decomposition}
Suppose that $|G|$ is invertible in~$R$, and let $M$ be an
$R\n$-functor.  Then
\[
\cK(M)(G)
=
\Delta_G(e_G^G)M(G),
\qquad
\cI(M)(G)
=
\Delta_G(1-e_G^G)M(G).
\]
In particular,
\[
M(G)
=
\cK(M)(G)\oplus\cI(M)(G)
\]
as $RB(G)$-modules.
\end{Prop}
\begin{proof}
By \Cref{lem:Burnside-action-compatibility} and
\Cref{rec:primitive-Burnside-idempotents},
$e_G^G\cdot M(G)\subseteq\cK(M)(G)$.
Conversely, if $m\in\cK(M)(G)$, then the formula for $1-e_G^G$ gives
$(1-e_G^G)\cdot m=0$, and hence
$\cK(M)(G)=e_G^G\cdot M(G)$.

Similarly, the formula for $1-e_G^G$ gives
$(1-e_G^G)\cdot M(G)\subseteq\cI(M)(G)$.
For $H<G$ and $n\in M(H)$,
\[
e_G^G\cdot{}\lind{G}_H(n)
=
{}\lind{G}_H
\bigl(({}\lres{H}_G(e_G^G))\cdot n\bigr)
=
0
\]
by \Cref{lem:Burnside-action-compatibility}.
Hence every properly induced element belongs to
$(1-e_G^G)\cdot M(G)$, so
$\cI(M)(G)=(1-e_G^G)\cdot M(G)$.

Since $RB(G)$ is commutative and $e_G^G$ is an idempotent, the asserted
direct-sum decomposition as $RB(G)$-modules follows.
\end{proof}
The following zero-detection criterion is the analogue of
\cite[Corollary~6.3]{Boltje98a}.
\begin{Cor}[Zero detection]
\label{cor:Burnside-zero-detection}
Suppose that $|G|$ is invertible in~$R$.  Let $M$ be an
$R\n$-functor and let $m\in M(G)$.  Then $m=0$ if and only if
$\Delta_G(e_G^G)m=0$ and ${}\lres{H}_G(m)=0$ for every strict
subgroup $H<G$.
\end{Cor}
\begin{proof}
This follows immediately from
\Cref{prop:Burnside-decomposition}.
\end{proof}

%
\section{The plus constructions}
\label{sec:plus-constructions}
%

This section develops the lower and upper plus constructions in the
change-of-rings framework of \Cref{sec:four-group-categories}. We
identify their module-theoretic realizations with the classical
cofixed and fixed descriptions, record their explicit actions and
multiplicative structures, and show that the corresponding
change-of-rings adjunctions restrict to Green functors.
%

%
\subsection{The lower plus construction}\
\label{ssec:lower-plus-construction}
%

We first define the lower plus construction at the level of modules.
\begin{Cons}[The lower-plus functor]
\label{con:lower-plus}
By \Cref{lem:four-group-categories}, we have an inclusion
$R\e\subseteq R\n$.  Moreover, the two locally unital algebras
have the same local units by \Cref{con:quiver-algebra}.

We define the lower plus functor
\[
\lp
:=
\iInd{R\n}{R\e}
\colon
\Fun(R\e)\longrightarrow\Fun(R\n).
\]
Thus, for an $R\e$-functor $A$,
$\lp(A)=R\n\otimes_{R\e}A$.
\end{Cons}
We next recall the classical cofixed construction and identify it
with the module-theoretic lower plus; compare
\cite[Sec.~2.2]{Boltje98a}.
\begin{Rec}[The cofixed lower-plus construction]
\label{rec:cofixed-lower-plus}
Let $A$ be an $R\e$-functor and let $G\in\mathfrak K$.
The group~$G$ acts on $\bigoplus\limits_{U\leq G}A(U)$ through
the isogation maps: for $g\in G$, the summand $A(U)$ is sent to
$A({}^gU)$ by ${}\liso{{}^gU}_U^{c_g}$.
	
We set
$A_+(G):=\Big(\bigoplus\limits_{U\leq G}A(U)\Big)_G$.
For $U\leq G$ and $x_U\in A(U)$, we write $[U,x_U]_G$ for the
class represented by~$x_U$ in the $U$-summand.
Thus
\[
[{}^gU,{}_{{}^gU}\iso_U^{c_g}(x_U)]_G
=
[U,x_U]_G
\qquad
(g\in G).
\]
Choosing one representative from each $G$-conjugacy class of
subgroups gives an $R$-module isomorphism
\[
A_+(G)
\cong
\bigoplus_{U\leq_G G}
A(U)_{N_G(U)}.
\]
Consequently, every element of $A_+(G)$ can be written as
\[
x
=
\sum_{U\leq_G G}[U,x_U]_G
\]
where $x_U\in A(U)$ and the class of $x_U$ in
$A(U)_{N_G(U)}$ is uniquely determined.
\end{Rec}
\begin{Prop}[Cofixed description of the lower plus]
\label{prop:lower-plus-cofixed}
Let $A$ be an $R\e$-functor and let $G\in\mathfrak K$.
There is an $R$-linear isomorphism
\[
s_G\colon
\lp(A)(G)
\isoto
A_+(G)
\]
determined by
$s_G({}\lind{G}_U\otimes x_U)=[U,x_U]_G$
for $U\leq G$ and $x_U\in A(U)$.
\end{Prop}
\begin{proof}
Define $t_G([U,x_U]_G):={}\lind{G}_U\otimes x_U$.
This is well defined, since for $g\in G$,
\[
{}\lind{G}_{{}^gU}\otimes
{}\liso{{}^gU}_U^{c_g}(x_U)
={}\lind{G}_{{}^gU}{}\!\liso{}_U^{c_g}\otimes x_U=
{}\lind{G}_U\otimes x_U
\]
by the balancing relation and
\Cref{not:basic-biset-morphisms}.

Conversely, let
$\eta\in\tm_{\mathcal N}(G,F)$.
Choose a factorization
$\eta={}\lind{G}_U\!e$
as in \Cref{lem:N-E-factorization}, and set
$s_G(\eta\otimes x_F):=[U,e(x_F)]_G$.

The uniqueness assertion of
\Cref{lem:N-E-factorization} and the cofixed relation of
\Cref{rec:cofixed-lower-plus} show that this is independent of
the chosen factorization. Extending $R$-linearly, it is enough to check the balancing
relation for a transitive morphism $r$ of $\mathcal E$. This follows from
the functoriality of $A$, since
\[
[U,(er)(x)]_G=[U,e(r(x))]_G.
\]
Finally,
$s_Gt_G=\id_{A_+(G)}$ and, if $\eta={}\lind{G}_U\!e$, then by the
balancing relation,
\[
t_Gs_G(\eta\otimes x)
=
{}\lind{G}_U\otimes e(x)
={}\lind{G}_U\!e\otimes x=
\eta\otimes x.
\]
Thus $s_G$ and $t_G$ are mutually inverse. Naturality in $A$ is immediate.
\end{proof}
\begin{Cons}
\label{con:lower-plus-action}
Let $A$ be an $R\e$-functor. For
$\eta\in R\n(F,G)$, define
\[
A_+(\eta)
:=
s_F\circ \lp(A)(\eta)\circ s_G^{-1}.
\]
Then the family $A_+(G)$ becomes an $R\n$-functor, and
$s\colon \lp(A)\isoto A_+$ is an $R\n$-isomorphism.

We henceforth identify $\lp(A)\cong A_+$ via~$s$.
If $f\colon A\to B$ is a morphism of $R\e$-functors, we write
$f_+\colon A_+\to B_+$ for the induced $R\n$-morphism. On generators,
\[
f_{+,G}([U,x]_G)=[U,f_U(x)]_G.
\]
\end{Cons}
\begin{Prop}[Explicit lower-plus action]
\label{prop:lower-plus-action}
Let $A$ be an $R\e$-functor.
\begin{enumerate}
\item[(a)]
Let $\alpha\colon H\to G$ be a morphism in $\mathcal I$,
let $V\leq H$, and let $y_V\in A(V)$.
Writing
$\ol{\alpha}_V\colon V\isoto\alpha(V)$
for the corestriction of $\alpha|_V$, we have
\[
\lind{G}_H^\alpha[V,y_V]_H
=
[\alpha(V),
{}\liso{\alpha(V)}_V^{\ol{\alpha}_V}(y_V)]_G.
\]

\item[(b)]
Let $\theta\colon H\to G$ be a morphism in $\mathcal R$,
let $U\leq G$, and let $x_U\in A(U)$.
For
$g\in\theta(H)\backslash G/U$, set
$K_g:=\theta^{-1}({}^gU)$.
Then
\[
\lres{H}_G^\theta[U,x_U]_G
=
\sum_{g\in\theta(H)\backslash G/U}
\left[
K_g,\,
{}\lres{K_g}_U^{\,c_{g^{-1}}\theta}(x_U)
\right]_H,
\]
where $c_{g^{-1}}\theta$ denotes the restriction to~$K_g$
of $c_{g^{-1}}\circ\theta$, whose image lies in~$U$.
\end{enumerate}
\end{Prop}
\begin{proof}
For~(a), under the isomorphism of
\Cref{prop:lower-plus-cofixed}, the element $[V,y_V]_H$
corresponds to $\lind{H}_V\otimes y_V$. By \Cref{rem:corestrictions},
$\liso{\alpha(V)}_V^{\ol{\alpha}_V}\in\e$.
Hence
\[
\lind{G}_H^\alpha
(\lind{H}_V\otimes y_V)
=
\lind{G}_{\alpha(V)}
\liso{}_V^{\ol{\alpha}_V}
\otimes y_V
=
\lind{G}_{\alpha(V)}
\otimes
\liso{\alpha(V)}_V^{\ol{\alpha}_V}(y_V),
\]
by \Cref{not:basic-biset-morphisms} and the balancing relation.
Applying $s_G$ proves~(a).

For~(b), under the same isomorphism,
$[U,x_U]_G$ corresponds to $\lind{G}_U\otimes x_U$.
By \Cref{rec:generalized-mackey-relation},
\[
\lres{H}_G^\theta\lind{}_U
=
\sum_{g\in\theta(H)\backslash G/U}
{}\lind{H}_{K_g}
\lres{}_U^{c_{g^{-1}}\theta},
\]
by \MS{1}--\MS{3},
$c_{g^{-1}}\theta\colon K_g\to U$
belongs to $\Mor(\mathcal R)$, so
$\lres{K_g}_U^{c_{g^{-1}}\theta}\in\e$.
Therefore, by the balancing relation,
\[
\lres{H}_G^\theta
(\lind{G}_U\otimes x_U)
=
\sum_{g\in\theta(H)\backslash G/U}
\lind{H}_{K_g}\otimes
\,{\lres{K_g\!}_U^{c_{g^{-1}}\theta}(x_U)}.
\]
Applying $s_H$ gives the formula in~(b).
\end{proof}
We next record the universal property of the lower plus construction.
\begin{Thm}[The lower-plus construction]
\label{thm:lower-plus-adjunction}
The lower plus functor
\[
(-)_{\scriptscriptstyle+}
\colon
\Fun(R\e)
\longrightarrow
\Fun(R\n)
\]
is left adjoint to
$\iRes{R\e}{R\n}$.
Thus, naturally in an $R\e$-functor $A$ and an
$R\n$-functor $M$, there is an $R$-linear isomorphism
\[
\Hom_{R\n}(\Alp,M)
\cong
\Hom_{R\e}
\bigl(
A,\iRes{R\e}{R\n}(M)
\bigr).
\]
\end{Thm}
\begin{proof}
By \Cref{con:lower-plus},
\Cref{Rec:change-of-rings-adjunctions} and
\Cref{rec:change-of-rings-locally-unital},
we have natural isomorphisms
\[
\Hom_{R\n}(\Alp,M)
\cong
\Hom_{R\n}
\bigl(
R\n\otimes_{R\e}A,M
\bigr)
\cong
\Hom_{R\e}
\bigl(
A,\iRes{R\e}{R\n}(M)
\bigr).
\]
Hence the lower plus functor is left adjoint to
$\iRes{R\e}{R\n}$.
\end{proof}
\begin{Rem}[Unit and counit of the lower-plus adjunction]
\label{rem:lower-plus-unit-counit}
The unit
\[
\etalp_A
\colon
A\longrightarrow
\iRes{R\e}{R\n}(\Alp)
\]
is, under the cofixed description of
\Cref{prop:lower-plus-cofixed}, given by
$\etalp_{A,G}(a)=[G,a]_G$.

The counit
\[
\epslp_M
\colon
\bigl(
\iRes{R\e}{R\n}(M)
\bigr)_{\scriptscriptstyle+}
\longrightarrow
M
\]
is given by
$\epslp_{M,G}([U,x_U]_G)=\lind{G}_U(x_U)$.
\end{Rem}
We now turn to the multiplicative structure. Assume that $A$ is a
Green $R\e$-functor.
\begin{Prop}[The lower-plus product]
\label{prop:lower-plus-product}
Let $A$ be a Green $R\e$-functor. For every
$G\in\mathfrak K$, there is a unital $R$-algebra structure on
$A_+(G)$ determined by
\[
[U,x_U]_G\cdot[V,y_V]_G
=
\sum_{g\in U\backslash G/V}
\left[
U^g\cap V,\,
\lres{U^g\cap V}_U^{c_g}(x_U)
\cdot
\lres{U^g\cap V}_V(y_V)
\right]_G,
\]
where $U,V\leq G$, $x_U\in A(U)$, and $y_V\in A(V)$.
The identity element is $[G,1_{A(G)}]_G$.
\end{Prop}
\begin{proof}
This is the lower-plus multiplication in the present setting.
The usual double-coset argument applies; compare
\cite[Theorem~7.3 and the paragraph following it]
{BoltjeRaggiCardenasValeroElizondo19}.
\end{proof}
\begin{Thm}[The Green lower plus construction]
\label{thm:lower-plus-Green}
Let $A$ be a Green $R\e$-functor. With the multiplication of
\Cref{prop:lower-plus-product}, the lower plus functor $A_+$ is a
Green $R\n$-functor.
	
Moreover, the lower plus construction restricts to a functor
\[
(-)_{\scriptscriptstyle+}
\colon
\grfun{R\e}
\longrightarrow
\grfun{R\n}.
\]
\end{Thm}
\begin{proof}
By \Cref{prop:lower-plus-product}, each $A_+(G)$ is a unital
$R$-algebra. By $R$-bilinearity, it is enough to verify the
multiplicativity of restriction and the Frobenius relations on
generators $[U,x_U]_G$. Substituting the formulas of
\Cref{prop:lower-plus-action,prop:lower-plus-product} and using
\Cref{rec:generalized-mackey-relation} gives the required identities.
Moreover, \Cref{prop:lower-plus-action}(b) gives
\[
\lres{I}_G^\theta[G,1_{A(G)}]_G
=
[I,1_{A(I)}]_I,
\]
so the restriction maps are unital. Hence $A_+$ is a Green
$R\n$-functor. Compare the corresponding calculation in
\cite[Theorem~7.3(a)]{BoltjeRaggiCardenasValeroElizondo19}.

Finally, let $f\colon A\to B$ be a morphism of Green
$R\e$-functors. By \Cref{con:lower-plus-action},
$f_{+,G}([U,x]_G)=[U,f_U(x)]_G$.
Since each $f_U$ is a unital algebra homomorphism commuting with
restrictions, \Cref{prop:lower-plus-product} shows that every
$f_{+,G}$ is a unital algebra homomorphism. Hence $f_+$ is a
morphism of Green $R\n$-functors.
\end{proof}
Having established the Green lower plus construction, we now record
its universal property.
\begin{Cor}[The Green lower-plus adjunction]
\label{cor:lower-plus-Green-adjunction}
The adjunction of \Cref{thm:lower-plus-adjunction} restricts to an
adjunction
\[
(-)_{\scriptscriptstyle+}
\colon
\grfun{R\e}
\rightleftarrows
\grfun{R\n}
\colon
\iRes{R\e}{R\n}.
\]
Equivalently, for a Green $R\e$-functor $A$ and a Green
$R\n$-functor $M$, there is a natural bijection
\[
\Hom_{R\n}^{\gr}(\Alp,M)
\cong
\Hom_{R\e}^{\gr}
\bigl(
A,\iRes{R\e}{R\n}(M)
\bigr).
\]
\end{Cor}
\begin{proof}
By \Cref{thm:lower-plus-Green,rem:Green-restrictions},
the two functors preserve Green functors.
By \Cref{def:Green-functor-morphisms}, it therefore remains to
check that the unit and counit of
\Cref{rem:lower-plus-unit-counit} are componentwise unital
$R$-algebra homomorphisms.

The unit has this property by
\Cref{prop:lower-plus-product}.
For the counit, it is enough to check multiplicativity on generators.
Let $x_U\in M(U)$ and $x_V\in M(V)$, and put
$K_g=U^g\cap V$. By \Cref{prop:lower-plus-product}, \Cref{rem:lower-plus-unit-counit},
the Frobenius relations, and the Mackey relation,
\[
\begin{aligned}
\epslp_{M,G}
\bigl([U,x_U]_G\cdot[V,x_V]_G\bigr)
&=
\sum_{g\in U\backslash G/V}
\lind{G}_{K_g}
\left(
\lres{K_g}_U^{c_g}(x_U)\cdot
\lres{K_g}_V(x_V)
\right) 
\\
&=
\lind{G}_V
\left(
\lres{V}_G\lind{}_U(x_U)\cdot x_V
\right) 
\\
&=
\lind{G}_U(x_U)\cdot\lind{G}_V(x_V) 
\\
&=\epslp_{M,G}
\bigl([U,x_U]_G\bigr) \cdot \epslp_{M,G}
\bigl([V,x_V]_G\bigr)
.
\end{aligned}
\]
Thus every $\epslp_{M,G}$ is multiplicative. Moreover,
$\epslp_{M,G}([G,1_{M(G)}]_G)=1_{M(G)}$, so it is unital.
\end{proof}

%
\subsection{The upper plus construction}\
\label{ssec:upper-plus-construction}
%

We now develop the upper plus construction in parallel with the lower
plus, beginning with its module-theoretic realization.
\begin{Cons}[The upper-plus functor]
\label{con:upper-plus}
By \Cref{lem:four-group-categories}, we have an inclusion
$R\w\subseteq R\n$.  Moreover,
\Cref{con:W-to-S-projection} gives a projection
$R\w\onto R\s$.

We define the upper plus functor
\[
\Up
:=
\iCo{R\n}{R\w}\!\iInf{}{R\s}
\colon
\Fun(R\s)\longrightarrow\Fun(R\n).
\]
Thus, for an $R\s$-functor $B$,
\[
\Up(B)
=
\iCo{R\n}{R\w}
\iInf{}{R\s}(B).
\]
Equivalently, by
\Cref{rec:change-of-rings-locally-unital},
\[
\Up(B)
=
\hm_{R\w}\!
\left(
R\n,\,
\lInf{R\w}_{R\s}(B)
\right)_{\mathrm{lu}}.
\]
\end{Cons}
This is the analogue of Coşkun's realization of the classical
upper plus construction as coinduction; see
\cite[Theorem~5.2]{Coskun07}.
To relate this construction to the classical fixed description, we
first identify its restriction to $R\e$.
\begin{Cor}[Restriction of the upper plus]
\label{cor:upper-plus-restriction-comparison}
For every $R\s$-functor $B$, there is a natural
$R\e$-isomorphism
\[\lRes{R\e}_{R\n}\Up(B)
\cong
\lCo{R\e}_{R\s}(B).
\]
\end{Cor}
\begin{proof}
This follows from
\Cref{prop:Ind-Coind-comparison,lem:Res-Inf-identity}.
\end{proof}

Although the upper plus construction is defined as an
$R\mathcal N$-functor, Corollary~\ref{cor:upper-plus-restriction-comparison}
identifies its restriction to $R\mathcal E$ with coinduction
from $R\mathcal S$. Consequently, the upper plus adjunction
naturally takes place between $R\mathcal E$-functors and
$R\mathcal S$-functors. The $R\mathcal N$-valued construction
provides the additional induction structure needed for the
mark morphism.

We next recall the classical fixed construction and identify it with
the module-theoretic upper plus; compare
\cite[Sec.~2.1]{Boltje98a}.
\begin{Rec}[The fixed upper-plus construction]
\label{rec:fixed-upper-plus}
Let $B$ be an $R\s$-functor and let $G\in\mathfrak K$.
As in \Cref{rec:cofixed-lower-plus}, the group~$G$ acts on
$\bigoplus\limits_{U\leq G}B(U)$ through isogations.

We set
\[
B^+(G)
:=
\left(
\bigoplus_{U\leq G}B(U)
\right)^G.
\]
Thus an element $\xi\in B^+(G)$ is a family
$\xi=(\xi^U)_{U\leq G}$ with $\xi^U\in B(U)$ satisfying
\begin{equation}
\label{eq:upper-plus-equivariance}
\xi^{{}^gU}
=
\liso{{}^gU}_U^{c_g}(\xi^U)
\qquad
(g\in G).
\end{equation}
Choosing one representative from each $G$-conjugacy class of
subgroups gives an $R$-module isomorphism
\[
B^+(G)
\cong
\bigoplus_{U\leq_GG}
B(U)^{N_G(U)}.
\]
For $\xi\in B(U)^{N_G(U)}$, we write $[U,\xi]^G$ for the element
of $B^+(G)$ corresponding to $\xi$ in the $U$-summand of the
above decomposition.
Hence every $\xi\in B^+(G)$ can be written uniquely as
\[
\xi
=
\sum_{U\leq_GG}[U,\xi^U]^G,
\qquad
\xi^U\in B(U)^{N_G(U)}.
\]
\end{Rec}
\begin{Lem}[Upper-plus evaluation]
\label{lem:upper-plus-evaluation}
Let $B$ be an $R\s$-functor and let $G\in\mathfrak K$.
Evaluation at the subgroup restrictions induces a natural
$R$-linear isomorphism
\[
\ev_G\colon
\hm_{R\s}
\bigl(
R\e\cdot\id_G,\,
B
\bigr)
\isoto
B^+(G)
\]
given by
$\ev_G(\epsilon)
=
\bigl(
\epsilon(\lres{U}_G)
\bigr)_{U\leq G}$.
\end{Lem}
\begin{proof}
Recall that by \Cref{con:quiver-algebra},
$R\e\cdot\id_G=\bigoplus_{H\in\mathfrak K}R\e(H,G)$,
the left $R\e$-module of morphisms with source $G$.
Let
$\epsilon\in\Hom_{R\s}(R\e\cdot\id_G,B)$.
For $g\in G$ and $U\leq G$, we have
\[
\epsilon(\lres{^gU}_G)
=
\epsilon\bigl(
\liso{^gU}_U^{c_g}\lres{}_G
\bigr)
=
\liso{^gU}_U^{c_g}
\bigl(\epsilon(\lres{U}_G)\bigr),
\]
by \Cref{not:basic-biset-morphisms} and $R\s$-linearity.
Hence $\ev_G(\epsilon)$ belongs to $B^+(G)$.

Conversely, let
$\xi=(\xi^U)_{U\leq G}\in B^+(G)$.
For a transitive morphism
$\lres{H}_G^\theta\in R\e\cdot\id_G$, we define
$\epsilon_\xi(\lres{H}_G^\theta)
:=
\lres{H}_{\theta(H)}^{\ol{\theta}}
\bigl(\xi^{\theta(H)}\bigr)$.
This is well defined. Indeed, if
$\lres{H}_G^\theta=\lres{H}_G^{\theta'}$, then, applying
\Cref{lem:right-free-transitive-uniqueness} to the presentations
determined by the corestrictions $\ol{\theta}$ and
$\ol{\theta'}$, there exist $h\in H$ and $g\in G$ such that
$\theta'(H)={}^g\theta(H)$ and
$\ol{\theta'}c_h=c_g\ol{\theta}$.
Using the fixed-family relation of
\Cref{rec:fixed-upper-plus} and
\Cref{not:basic-biset-morphisms}, we obtain
\[
{}\lres{H}_{\theta'(H)}^{\ol{\theta'}}
\bigl(\xi^{\theta'(H)}\bigr)
=
\lres{H}_{\theta(H)}^{\ol{\theta}}
\bigl(\xi^{\theta(H)}\bigr).
\]
Since the transitive morphisms form an $R$-basis of
$R\e\cdot\id_G$ by \Cref{not:transitive-morphisms},
$\epsilon_\xi$ extends uniquely $R$-linearly to
$R\e\cdot\id_G$.

It is $R\s$-linear. By $R$-linearity, it is enough to check this
for a surjective $\sigma\colon K\onto H$ in $\Mor(\mathcal R)$ and a
transitive morphism $\lres{H}_G^\theta$. Since $\sigma$ is surjective,
$(\theta\sigma)(K)=\theta(H)$, and the corestriction of
$\theta\sigma$ is $\ol{\theta}\sigma$. Hence, by
\Cref{not:basic-biset-morphisms},
\[
\epsilon_\xi\!
\bigl(
\lres{K}_H^\sigma\!\lres{}_G^\theta
\bigr)\!
=
\epsilon_\xi\!
\bigl(
\lres{K}_G^{\theta\sigma}
\bigr) \!
=\!
\lres{K}_{\theta(H)}^{\ol{\theta}\sigma}
\bigl(\xi^{\theta(H)}\bigr) \!
=\!
\lres{K}_H^\sigma\!
\lres{}_{\theta(H)}^{\ol{\theta}}
\bigl(\xi^{\theta(H)}\bigr)\!
=\!
\lres{K}_H^\sigma\!
\epsilon_\xi(\lres{H}_G^\theta).
\]
Thus $\epsilon_\xi$ is $R\s$-linear.

Finally, for $U\leq G$, the corestriction of the inclusion
$U\hookrightarrow G$ is $\id_U$, so
$\epsilon_\xi(\lres{U}_G)=\xi^U$.
Hence $\ev_G(\epsilon_\xi)=\xi$.
Conversely, $R\s$-linearity gives
\[
\epsilon(\lres{H}_G^\theta)
=
\lres{H}_{\theta(H)}^{\ol{\theta}}
\bigl(\epsilon(\lres{\theta(H)}_G)\bigr),
\]
so $\epsilon_{\ev_G(\epsilon)}=\epsilon$.
Therefore $\ev_G$ is an isomorphism.
Naturality in $B$ is immediate. 
\end{proof}
\begin{Prop}[Fixed-point description of the upper plus]
\label{prop:upper-plus-fixed}
Let $B$ be an $R\s$-functor and let $G\in\mathfrak K$.
There is a natural $R$-linear isomorphism
\[
t_G\colon
\Up(B)(G)
\isoto
B^+(G)
\]
given by
\[
t_G(x)
=
\sum_{U\leq_GG}
[U,x(\lres{U}_G)]^G.
\]
\end{Prop}
\begin{proof}
By \Cref{cor:upper-plus-restriction-comparison},
the component $\Up(B)(G)$ is naturally identified with
\[
\hm_{R\s}
\bigl(
R\e\cdot\id_G,\,
B
\bigr),
\]
where the identification is given by restriction as in
\Cref{prop:Ind-Coind-comparison}.

Now apply
\Cref{lem:upper-plus-evaluation}.
It gives an $R$-linear isomorphism to $B^+(G)$, and its
evaluation formula sends $x$ to the invariant family whose
$U$-coordinate is $x(\lres{U}_G)$.
Using the notation of
\Cref{rec:fixed-upper-plus}, this is precisely
\[
\sum_{U\leq_GG}[U,x(\lres{U}_G)]^G.
\]
Naturality follows from the naturality of
\Cref{cor:upper-plus-restriction-comparison} and
\Cref{lem:upper-plus-evaluation}.
\end{proof}
\begin{Cons}
\label{con:upper-plus-action}
Let $B$ be an $R\s$-functor. For $\eta\in R\n(F,G)$, define
\[
B^+(\eta)
:=
t_F\circ \Up(B)(\eta)\circ t_G^{-1}.
\]
Then the family $B^+(G)$ becomes an $R\n$-functor, and
$t\colon \Up(B)\isoto B^+$ is an $R\n$-isomorphism.

We henceforth write $B^+$ for this $R\n$-functor and identify
$\Up(B)\cong B^+$ via~$t$.
If $f\colon B\to B'$ is a morphism of $R\s$-functors, we write
$f^+\colon B^+\to (B')^+$ for the induced $R\n$-morphism. On
coordinates,
\[
\bigl(f_G^+(\xi)\bigr)^U=f_U(\xi^U).
\]
\end{Cons}
\begin{Prop}[Explicit upper-plus action]
\label{prop:upper-plus-action}
Let $B$ be an $R\s$-functor.

\begin{enumerate}
\item[(a)]
Let $\theta\colon H\to G$ be a morphism in $\mathcal R$
and let $\xi\in B^+(G)$.
For $V\leq H$, let
$\ol{\theta}_V\colon V\onto\theta(V)$ be the corestriction
of $\theta|_V$. Then
\[
\bigl(\lres{H}_G^\theta(\xi)\bigr)^V
=
\lres{V}_{\theta(V)}^{\ol{\theta}_V}
\bigl(\xi^{\theta(V)}\bigr).
\]
\item[(b)]
Let $\alpha\colon I\to F$ be a morphism in $\mathcal I$
and let $\eta\in B^+(I)$.
For $U\leq F$ and
$g\in U\backslash F/\alpha(I)$, set
$K_g
:=
\alpha^{-1}({U^g})$ and
$\gamma_g
:=
c_g\alpha|_{K_g}\colon K_g\to U$.
Then
\[
\bigl(\lind{F}_I^\alpha(\eta)\bigr)^U
=
\sum_{\substack{
		g\in U\backslash F/\alpha(I)\\
		U\leq{}^g\alpha(I)
}}
\liso{U}_{K_g}^{\gamma_g}
\bigl(\eta^{K_g}\bigr).
\]
For every term occurring in the sum, $\gamma_g$ is an
isomorphism.
\end{enumerate}
\end{Prop}
\begin{proof}
Under the identification of
\Cref{prop:upper-plus-fixed}, we regard elements of $B^+(G)$ as
$R\w$-homomorphisms on $R\n\cdot\id_G$.
The $R\n$-action is given by right multiplication in the argument.

For~(a), if $V\leq H$, then
\[
\bigl(\lres{H}_G^\theta(\xi)\bigr)^V
=
\xi\bigl(\lres{V}_H\lres{}_G^\theta\bigr)
=
\xi\bigl(
\lres{V}_{\theta(V)}^{\ol{\theta}_V}
\lres{}_G
\bigr)
=
\lres{V}_{\theta(V)}^{\ol{\theta}_V}
\bigl(\xi^{\theta(V)}\bigr),
\]
where the last equality uses $R\w$-linearity, since
$\ol{\theta}_V$ is surjective and hence
$\lres{V}_{\theta(V)}^{\ol{\theta}_V}\in\mathcal S\subseteq\mathcal W$.
This proves~(a).

Now we prove~(b).
For $g\in U\backslash F/\alpha(I)$, put
$K_g=\alpha^{-1}(U^g)$ and
$\gamma_g=c_g\alpha|_{K_g}\colon K_g\to U$.
Since $K_g\hookrightarrow I$ and conjugation by $g$ belong to
$\Mor(\mathcal I)$ by \MS{1} and \MS{2},
\Cref{rem:corestrictions} and \MS{1} give
$\gamma_g\in\Mor(\mathcal I)$.
Then
\[
\begin{aligned}
\bigl(\lind{F}_I^\alpha(\eta)\bigr)^U
&=
\eta\bigl(\lres{U}_F\lind{}_I^\alpha\bigr)
&&\text{by the coinduction action}
\\
&=
\sum_{g\in U\backslash F/\alpha(I)}
\eta\bigl(
\lind{U}_{K_g}^{\gamma_g}\!\lres{}_I
\bigr)
&&\text{by \Cref{rec:generalized-mackey-relation}}
\\
&=
\sum_{g\in U\backslash F/\alpha(I)}
\lind{U}_{K_g}^{\gamma_g}
\bigl(\eta(\lres{K_g}_I)\bigr)
&&\text{by $R\w$-linearity of $\eta$}
\\
&=
\sum_{g\in U\backslash F/\alpha(I)}
\lind{U}_{K_g}^{\gamma_g}
\bigl(\eta^{K_g}\bigr)
&&\text{by \Cref{prop:upper-plus-fixed}.}
\end{aligned}
\]
Now
$\gamma_g(K_g)=U\cap{}^g\alpha(I)$.
If $U\not\leq{}^g\alpha(I)$, then
$\gamma_g(K_g)<U$, so
$\lind{U}_{K_g}^{\gamma_g}\in\Delta_\w$ by
\Cref{lem:Delta-W-transitive}.
Since $R\Delta_\w$ acts trivially after inflation by
\Cref{not:W-restriction-inflation}, the corresponding term vanishes.

If $U\leq{}^g\alpha(I)$, then
$\gamma_g(K_g)=U$.
Since $\gamma_g$ is injective, it is an isomorphism, and
\[
\lind{U}_{K_g}^{\gamma_g}
=
\liso{U}_{K_g}^{\gamma_g}
\]
by \Cref{not:basic-biset-morphisms}.
Therefore
\[
\bigl(\lind{F}_I^\alpha(\eta)\bigr)^U
=
\sum_{\substack{
	g\in U\backslash F/\alpha(I)\\
	U\leq{}^g\alpha(I)
}}
\liso{U}_{K_g}^{\gamma_g}
\bigl(\eta^{K_g}\bigr),
\]
which proves~(b).
\end{proof}
We next record the universal property of the upper plus construction.
\begin{Thm}[The upper-plus adjunction]
\label{thm:upper-plus-adjunction}
The functor
\[
\iRes{R\e}{R\n} \circ \Up
\colon
\Fun(R\s)
\longrightarrow
\Fun(R\e)
\]
is right adjoint to
$\iRes{R\s}{R\e}$.
Thus, naturally in an $R\e$-functor $A$ and an
$R\s$-functor $B$, there is a natural $R$-linear isomorphism
\[
\Hom_{R\s}
\bigl(
\iRes{R\s}{R\e}(A),B
\bigr)
\cong
\Hom_{R\e}
\bigl(
A,\iRes{R\e}{R\n}(\Bup)
\bigr).
\]
\end{Thm}
\begin{proof}
By \Cref{cor:upper-plus-restriction-comparison} and the
restriction--coinduction adjunction of
\Cref{Rec:change-of-rings-adjunctions,rec:change-of-rings-locally-unital},
\[
\hm_{R\s}
\bigl(
\iRes{R\s}{R\e}(A),B
\bigr)
\cong
\hm_{R\e}
\bigl(
A,{}\lCo{R\e}_{R\s}(B)
\bigr)
\cong
\hm_{R\e}
\bigl(
A, \iRes{R\e}{R\n}(B^+)
\bigr).
\]
\end{proof}
\begin{Rem}[Unit and counit of the upper-plus adjunction]
\label{rem:upper-plus-unit-counit}
The unit is
\[
\etaup_A
\colon
A
\longrightarrow
\iRes{R\e}{R\n}
\Bigl(
\Bup[\big(\iRes{R\s}{R\e}(A)\big)]
\Bigr),
\]
under the fixed description of
\Cref{prop:upper-plus-fixed}, given by
$\bigl(\etaup_{A,G}(a)\bigr)^U
=
\lres{U}_G(a)$.

The counit
\[
\epsup_B
\colon
\iRes{R\s}{R\n}
(\Bup)
\longrightarrow
B
\]
is the projection onto the top coordinate:
$\epsup_{B,G}(\xi)=\xi^G$.
\end{Rem}
We now turn to the multiplicative structure. Assume that $B$ is a
Green $R\s$-functor.
\begin{Cons}[The upper-plus product]
\label{con:upper-plus-product}
Let $B$ be a Green $R\s$-functor.
For $\xi,\eta\in B^+(G)$, define
\[
(\xi\eta)^U
:=
\xi^U\cdot\eta^U
\qquad
(U\leq G).
\]
The identity element is the family $1_{B^+(G)}$ given by
$(1_{B^+(G)})^U
=
1_{B(U)}$ for every
$U\leq G$.

The product and identity defined above satisfy the equivariance
condition of Eq. (\ref{eq:upper-plus-equivariance}), since isogations
are unital algebra homomorphisms.
Associativity and the unit identities hold coordinatewise.
Thus $B^+(G)$ is a unital $R$-algebra.
\end{Cons}
\begin{Thm}[The Green upper-plus construction]
\label{thm:upper-plus-Green}
Let $B$ be a Green $R\s$-functor. With the multiplication of
\Cref{con:upper-plus-product},
$\Bup$
is a Green $R\n$-functor.

Moreover,
\[
(-)^{\scriptscriptstyle+}
\colon
\grfun{R\s}
\longrightarrow
\grfun{R\n}
\]
is a functor.
\end{Thm}
\begin{proof}
By \Cref{con:upper-plus-product}, each $B^+(G)$ is a unital
$R$-algebra. The formulas of \Cref{prop:upper-plus-action},
together with the coordinatewise multiplication, give the
restriction and Frobenius axioms. Hence $B^+$ is a Green
$R\n$-functor; compare
\cite[Theorem~7.3(b)]{BoltjeRaggiCardenasValeroElizondo19}.

Finally, let $f\colon B\to B'$ be a morphism of Green
$R\s$-functors. By \Cref{con:upper-plus-action},
\[
\bigl(f_G^+(\xi)\bigr)^U=f_U(\xi^U).
\]
Since each $f_U$ is a unital $R$-algebra homomorphism,
\Cref{con:upper-plus-product} shows that every $f_G^+$ is a
unital $R$-algebra homomorphism. Hence $f^+$ is a morphism of
Green $R\n$-functors.
\end{proof}
Having established the Green upper plus construction, we now record
its universal property.
\begin{Cor}[The Green upper-plus adjunction]
\label{cor:upper-plus-Green-adjunction}
The adjunction of \Cref{thm:upper-plus-adjunction} restricts to an
adjunction
\[
\iRes{R\s}{R\e}
\colon
\grfun{R\e}
\rightleftarrows
\grfun{R\s}
\colon
\iRes{R\e}{R\n}\circ(-)^{\scriptscriptstyle+}.
\]
Equivalently, for a Green $R\e$-functor $A$ and a Green
$R\s$-functor $B$, there is a natural bijection
\[
\Hom_{R\s}^{\gr}
\bigl(
\iRes{R\s}{R\e}(A),B
\bigr)
\cong
\Hom_{R\e}^{\gr}
\bigl(
A,\iRes{R\e}{R\n}(\Bup)
\bigr).
\]
\end{Cor}
\begin{proof}
By
\Cref{thm:upper-plus-Green,rem:Green-restrictions},
the two functors preserve Green functors.
By \Cref{def:Green-functor-morphisms}, it remains to check that the
unit and counit of \Cref{rem:upper-plus-unit-counit} are
componentwise unital $R$-algebra homomorphisms.

The components of the unit are built from restriction maps, hence
are unital algebra homomorphisms. The counit is projection onto the
top coordinate, which is a unital algebra homomorphism for the
coordinatewise multiplication of
\Cref{con:upper-plus-product}.
\end{proof}

%
%
\section{The mark map and its inverse}
\label{sec:mark-map}
%

We now compare the lower and upper plus constructions.
The adjunctions of
\Cref{thm:lower-plus-adjunction,thm:upper-plus-adjunction}
yield the mark map $A_+\to A^+$.
We derive its explicit formula and, when the relevant group
orders are invertible in $R$, compute its inverse by
Möbius inversion.

%
\subsection{Construction and basic properties}
\label{ssec:mark-map}
%

\begin{Not}
\label{not:upper-plus-E-functor}
For an $R\e$-functor $A$, we write
\[
A^+
:=
\Bup[(\iRes{R\s}{R\e}(A))].
\]
Thus $A^+$ denotes the upper plus construction applied to the
underlying $R\s$-functor of~$A$.
\end{Not}
\begin{Cons}[The Brauer projection]
\label{con:Brauer-projection}
Let $A$ be an $R\e$-functor and let $G\in\mathfrak K$.
We define an $R$-linear map
$\br^A
\colon
\Alp
\longrightarrow
A$
by
\[
\br_G^A([U,x_U]_G)
=
\begin{cases}
	x_U,& U=G,\\
	0,& U<G.
\end{cases}
\]
Equivalently, under the decomposition of
\Cref{rec:cofixed-lower-plus}, the map $\br_G^A$ is the
projection onto the summand
$A(G)_{N_G(G)}=A(G)$.
\end{Cons}
By \Cref{prop:lower-plus-action}, these projections are compatible
with restriction along the surjective morphisms of $\mathcal R$ and
therefore assemble into an $R\s$-morphism
\[
\br^A\colon
\lRes{R\s}_{R\n}(A_+)
\longrightarrow
\lRes{R\s}_{R\e}(A).
\]
\begin{Prop}[The plus adjunction correspondence]
\label{prop:plus-adjunction-correspondence}
Let $A$ be an $R\e$-functor.
The two plus adjunctions induce an $R$-linear isomorphism
\[
\Omega^A\colon
\Hom_{R\s}
\bigl(
\iRes{R\s}{R\e}(A),
\iRes{R\s}{R\e}(A)
\bigr)
\isoto
\Hom_{R\n}(\Alp,A^+).
\]
If $\phi$ corresponds to
$\psi=\Omega^A(\phi)$, then, for every
$G\in\mathfrak K$, $U\leq G$, and $x\in A_+(G)$,
\[
\bigl(\psi_G(x)\bigr)^U
=
\phi_U
\left(
\br_U^A({}\lres{U}_G(x))
\right).
\]
If $A$ is a Green $R\e$-functor, this correspondence
restricts to a bijection
\[
\Hom_{R\s}^{\gr}
\bigl(
\iRes{R\s}{R\e}(A),
\iRes{R\s}{R\e}(A)
\bigr)
\cong
\Hom_{R\n}^{\gr}(A_+,A^+).
\]
\end{Prop}

\begin{proof}
We get
\[
\hm_{R\s}
\bigl(
\iRes{R\s}{R\e}(A),
\iRes{R\s}{R\e}(A)
\bigr)
\congby{\ref{thm:upper-plus-adjunction}}
\hm_{R\e}
\bigl(
A,
\iRes{R\e}{R\n}(A^+)
\bigr)\congby{\ref{thm:lower-plus-adjunction}}\hm_{R\n}(A_+,A^+).
\]
The composite is the asserted isomorphism $\Omega^A$.

We now establish the coordinate formula.
Let
$\widetilde\phi\colon A\to\iRes{R\e}{R\n}(A^+)$
correspond to $\phi$ under the upper-plus adjunction. By
\Cref{rem:upper-plus-unit-counit}, $\bigl(\widetilde\phi_G(a)\bigr)^U
=
\phi_U\bigl(\lres{U}_G(a)\bigr)$,
so in particular
$\bigl(\widetilde\phi_U(a)\bigr)^U=\phi_U(a)$.

Let $\psi=\Omega^A(\phi)$. By
\Cref{rem:lower-plus-unit-counit},
$\psi_U([V,a]_U)
=
\lind{U}_V\bigl(\widetilde\phi_V(a)\bigr)$.
Thus, by \Cref{prop:upper-plus-action},
\[
\bigl(\psi_U([V,a]_U)\bigr)^U
=
\begin{cases}
	\phi_U(a),&V=U,\\
	0,&V<U,
\end{cases}
\]
since $U\leq{}^gV$ for some $g\in U$ only when $V=U$.
Hence, by $R$-linearity and the definition of $\br_U^A$,
\begin{equation}
	\label{eq:plus-adjunction-top-coordinate}
	\bigl(\psi_U(y)\bigr)^U
	=
	\phi_U\bigl(\br_U^A(y)\bigr)
	\qquad
	(y\in A_+(U)).
\end{equation}
Now let $x\in A_+(G)$ and $U\leq G$. Since $\psi$ is an
$R\n$-morphism,
\[
\lres{U}_G\bigl(\psi_G(x)\bigr)
=
\psi_U\bigl(\lres{U}_G(x)\bigr).
\]
Taking the $U$-coordinate and using
\Cref{prop:upper-plus-action,eq:plus-adjunction-top-coordinate}
give
\[
\bigl(\psi_G(x)\bigr)^U
=
\phi_U\left(
\br_U^A\bigl(\lres{U}_G(x)\bigr)
\right),
\]
as claimed.

If $A$ is Green, the last assertion follows from
\Cref{cor:lower-plus-Green-adjunction} and	
\Cref{cor:upper-plus-Green-adjunction}, since these are
restrictions of the corresponding adjunctions for
ordinary functors.
\end{proof}
\begin{Cons}[The mark map]
\label{con:mark-map}
Let $A$ be an $R\e$-functor.
Following the module-theoretic interpretation of
\cite[Corollary~5.10]{Coskun07}, we define the \emph{mark map}
\[
\mk^A\colon A_+\longrightarrow A^+
\]
to be the $R\n$-morphism corresponding under
\Cref{prop:plus-adjunction-correspondence} to the identity
morphism of $\iRes{R\s}{R\e}(A)$.  Thus
\[
\mk^A
=
\Omega^A
\bigl(
\id_{(\iRes{R\s}{R\e}(A))}
\bigr).
\]
\end{Cons}
\begin{Prop}[The mark formula]
\label{prop:mark-formula}
Let $A$ be an $R\e$-functor, let $G\in\mathfrak K$, and let
$x\in A_+(G)$.
For every $U\leq G$,
\[
\bigl(\mk_G^A(x)\bigr)^U
=
\br_U^A(\lres{U}_G(x)).
\]
In particular, for $V\leq G$ and $x_V\in A(V)$,
\[
\bigl(\mk_G^A([V,x_V]_G)\bigr)^U
=
\sum_{\substack{
		gV\in G/V\\
		U\leq{}^gV
}}
\lres{U}_{{}^gV}
{}\liso{}_V^{c_g}(x_V).
\]
\end{Prop}
\begin{proof}
The first formula follows from
\Cref{prop:plus-adjunction-correspondence} applied to
$\phi=\id_{(\iRes{R\s}{R\e}(A))}$, together with
\Cref{con:mark-map}.
	
For the second formula, apply
\Cref{prop:lower-plus-action} to the inclusion
$U\hookrightarrow G$.
We obtain
\[
\lres{U}_G[V,x_V]_G
=
\sum_{g\in U\backslash G/V}
\left[
U\cap{}^gV,\,\lres{U\cap{}^gV}_V^{c_{g^{-1}}}(x_V)
\right]_U.
\]
By \Cref{con:Brauer-projection}, only summands
with $U\cap{}^gV=U$, equivalently
$U\leq{}^gV$, can contribute to $\br_U^A$.
For such a~$g$,
\Cref{not:basic-biset-morphisms} gives
\[
\lres{U}_V^{c_{g^{-1}}}
=
\lres{U}_{{}^gV}\liso{}_V^{c_g}.
\]
Moreover, if $U\leq{}^gV$, then
$UgV=gV$.
Hence the nonzero double-coset terms are in bijection with the cosets
$gV\in G/V$ satisfying $U\leq{}^gV$.
Substituting into the first formula proves the result.
\end{proof}
\begin{Rem}[Comparison with the classical mark map]
\label{rem:classical-mark-map}
In the classical subgroup setting, the formula of
\Cref{prop:mark-formula} is Boltje's mark homomorphism.
In particular, the maps $\br_U^A$ and $\mk_G^A$ correspond to
the maps denoted $\pi_U^A$ and $\rho_G^A$, respectively, in
\cite[Sec.~2.3]{Boltje98a}.
The module-theoretic characterization of
\Cref{con:mark-map} agrees with
\cite[Corollary~5.10]{Coskun07}.
\end{Rem}
\begin{Thm}[The Green mark map]
\label{thm:Green-mark-map}
Let $A$ be a Green $R\e$-functor.
Then
\[
\mk^A\colon A_+\longrightarrow A^+
\]
is a morphism of Green $R\n$-functors.
Moreover, the maps $\mk^A$ are natural in~$A$.
\end{Thm}
\begin{proof}
By \Cref{prop:plus-adjunction-correspondence},
the correspondence $\Omega^A$ restricts to Green functors.
Since the identity of
$\iRes{R\s}{R\e}(A)$ is a morphism of Green
$R\s$-functors, its image $\mk^A$ is a morphism
of Green $R\n$-functors by \Cref{con:mark-map}.

Naturality in $A$ follows from the coordinate formula of
\Cref{prop:mark-formula} and the naturality of the
Brauer projections.
\end{proof}
%

%
\subsection{The inverse mark formula}\
\label{ssec:inverse-mark-formula}
%

We now study the invertibility of the mark map, following
Boltje~\cite[Proposition~2.4]{Boltje98a}.
We use the subgroup Möbius function $\mu$ on the subgroup
lattice of a finite group $G$, as recalled in
\Cref{rec:primitive-Burnside-idempotents}.
\begin{Cons}[The inverse-mark operator]
\label{con:inverse-mark-operator}
Let $A$ be an $R\e$-functor and let $G\in\mathfrak K$.
We define an $R$-linear map
$\km_G^A\colon A^+(G)\longrightarrow A_+(G)$
by
\[
\km_G^A(\xi)
:=
\sum_{U\leq V\leq G}
|U|\mu(U,V)
\left[
U,\,
\lres{U}_V(\xi^V)
\right]_G.
\]
Here the sum runs over all pairs of subgroups $U\leq V$ of~$G$,
not over conjugacy-class representatives.
\end{Cons}
\begin{Lem}[Möbius inversion for the mark map]
\label{lem:mark-Mobius-inversion}
Let $A$ be an $R\e$-functor and let $G\in\mathfrak K$.
Then
$\km_G^A\mk_G^A=|G|\,\id_{A_+(G)}$
and
$\mk_G^A\km_G^A=
|G|\,\id_{A^+(G)}$.
\end{Lem}
\begin{proof}
For fixed $G$, the calculation uses only the subgroup lattice of
$G$, the mark formula of \Cref{prop:mark-formula}, and the cofixed
and fixed relations of
\Cref{rec:cofixed-lower-plus,rec:fixed-upper-plus}.
Thus the Möbius-inversion argument of
\cite[Proposition~2.4]{Boltje98a} applies verbatim and gives both
identities.
\end{proof}
\begin{Prop}[The inverse mark formula]
\label{prop:inverse-mark-formula}
Let $A$ be an $R\e$-functor and let $G\in\mathfrak K$.
Suppose that $|G|$ is invertible in~$R$.
Then
\[
\mk_G^A\colon A_+(G)\isoto A^+(G)
\]
is an $R$-linear isomorphism.
For every $\xi\in A^+(G)$,
\begin{equation}
\label{eq:mark-inverse}
(\mk_G^A)^{-1}(\xi)
=
\frac{1}{|G|}
\sum_{U\leq V\leq G}
|U|\mu(U,V)
\left[
U,\,
{}\lres{U}_V(\xi^V)
\right]_G.
\end{equation}
\end{Prop}
\begin{proof}
By \Cref{lem:mark-Mobius-inversion},
\[
\km_G^A\mk_G^A
=
|G|\id_{A_+(G)}
\qquad\text{and}\qquad
\mk_G^A\km_G^A
=
|G|\id_{A^+(G)}.
\]
Since $|G|$ is invertible in~$R$, the map
$|G|^{-1}\km_G^A$ is both a left and a right inverse of
$\mk_G^A$.
The displayed formula follows from
\Cref{con:inverse-mark-operator}.
\end{proof}
\begin{Cor}[Invertibility of the mark map]
\label{cor:mark-isomorphism}
Suppose that $|G|$ is invertible in $R$ for every
$G\in\mathfrak K$.
Then, for every $R\e$-functor $A$, the mark map
\[
\mk^A\colon A_+\isoto A^+
\]
is an isomorphism of $R\n$-functors.

If $A$ is a Green $R\e$-functor, then $\mk^A$ is an
isomorphism of Green $R\n$-functors.
\end{Cor}
\begin{proof}
By \Cref{prop:inverse-mark-formula}, every component
$\mk_G^A$ is bijective.
Since $\mk^A$ is an $R\n$-morphism by
\Cref{con:mark-map}, its componentwise inverse is
an $R\n$-morphism.
Hence $\mk^A$ is an isomorphism of $R\n$-functors.

If $A$ is Green, then $\mk^A$ is a morphism of Green
$R\n$-functors by \Cref{thm:Green-mark-map}.
Each component is a unital algebra isomorphism, so its inverse
is also a unital algebra homomorphism.
Thus the inverse is a morphism of Green $R\n$-functors by
\Cref{def:Green-functor-morphisms}.
\end{proof}

%
\section{Canonical induction}
\label{sec:canonical-induction}
%
We now study canonical induction in the setting of Green
categories. We first define canonical induction as a section
of the linearization map. We then use the upper-plus adjunction
and the inverse mark map to construct a candidate for such a
section, and determine when it is a canonical induction formula.
We also study its multiplicativity.

We begin with the morphism from the lower plus construction to an
$R\n$-functor determined by the induction--restriction adjunction.
\begin{Cons}[The linearization map]
\label{con:linearization-map}
Let $A$ be an $R\e$-functor and let $M$ be an
$R\n$-functor.
By \Cref{thm:lower-plus-adjunction}, there is a
natural $R$-linear isomorphism
\[
\Hom_{R\n}(A_+,M)
\cong
\Hom_{R\e}
\bigl(
A,\iRes{R\e}{R\n}(M)
\bigr).
\]

For an $R\e$-morphism
$\nu\colon A\to\iRes{R\e}{R\n}(M)$,
we denote the corresponding $R\n$-morphism by
\[
\lin^\nu\colon A_+\longrightarrow M
\]
and call it the \emph{linearization map associated with $\nu$}.

Under the cofixed description of
\Cref{prop:lower-plus-cofixed}, it is given by
\[
\lin_G^\nu([U,x_U]_G)
=
\lind{G}_U\bigl(\nu_U(x_U)\bigr)
\]
for $U\leq G$ and $x_U\in A(U)$.

If $A$ and $M$ are Green functors and $\nu$ is a
Green $R\e$-morphism, then $\lin^\nu$ is a
Green $R\n$-morphism by
\Cref{cor:lower-plus-Green-adjunction}.

In Boltje's classical setting, when $A$ is a restriction
subfunctor of $M$ and $\nu$ is the inclusion, $\lin^\nu$
is the induction morphism $b^{M,A}$ of
\cite[Sec.~3.1]{Boltje98a}.
\end{Cons}
\begin{Def}[Canonical induction formula]
\label{def:canonical-induction-formula}
Let $A$ be an $R\e$-functor, let $M$ be an
$R\n$-functor, and let
$\nu\colon A\to\iRes{R\e}{R\n}(M)$
be an $R\e$-morphism.

A \emph{canonical induction formula for $\lin^\nu$}
is an $R\e$-morphism
\[
\alpha\colon
\iRes{R\e}{R\n}(M)
\longrightarrow
\iRes{R\e}{R\n}(A_+)
\]
such that
\[
\lin^\nu\circ\alpha
=
\id_{\iRes{R\e}{R\n}(M)}.
\]
\end{Def}
We next use the restriction--coinduction adjunction to construct a
morphism into the upper plus construction.
\begin{Cons}[The tom Dieck map]
\label{con:tom-Dieck-map}
Let $M$ be an $R\n$-functor and let $A$ be an
$R\e$-functor. By \Cref{thm:upper-plus-adjunction}, there is a 
natural $R$-linear isomorphism
\[
\Hom_{R\s}
\bigl(
{}\lRes{R\s}_{R\n}(M),
{}\lRes{R\s}_{R\e}(A)
\bigr)
\cong
\Hom_{R\e}
\bigl(
{}\lRes{R\e}_{R\n}(M),
{}\lRes{R\e}_{R\n}(A^+)
\bigr).
\]
Under this isomorphism, given an $R\s$-morphism
$\pi\colon
{}\lRes{R\s}_{R\n}(M)
\longrightarrow
{}\lRes{R\s}_{R\e}(A)$,
there exists a corresponding $R\e$-morphism
$\die^\pi\colon
{}\lRes{R\e}_{R\n}(M)
\longrightarrow
{}\lRes{R\e}_{R\n}(A^+)$
which we call the \emph{tom Dieck map associated with $\pi$}.

Under the fixed description of
\Cref{prop:upper-plus-fixed}, its coordinates are
\[
\bigl(\die_G^\pi(m)\bigr)^U
=
\pi_U\bigl({}\lres{U}_G(m)\bigr)
\]
for $U\leq G$. Equivalently,
\[
\die_G^\pi(m)
=
\sum_{U\leq_GG}
\left[
U,\,
\pi_U\bigl({}\lres{U}_G(m)\bigr)
\right]^G.
\]
Indeed, the top coordinate is given by
\Cref{rem:upper-plus-unit-counit}, and the remaining coordinates
follow from restriction compatibility and
\Cref{prop:upper-plus-action}.
\end{Cons}
We now use the inverse mark map to construct a candidate
for a canonical induction formula. By
\Cref{cor:mark-isomorphism}, the mark map is invertible
whenever $|G|$ is invertible in $R$ for every
$G\in\mathfrak K$.
\begin{Cons}[The canonical induction map]
\label{con:canonical-induction-map}
Let $M$ be an $R\n$-functor and let $A$ be an
$R\e$-functor. Suppose that
\[
\mk^A\colon A_+\longrightarrow A^+
\]
is an isomorphism of $R\n$-functors, and let
$\pi\colon
{}\lRes{R\s}_{R\n}(M)
\longrightarrow
{}\lRes{R\s}_{R\e}(A)$
be an $R\s$-morphism.
We define the \emph{canonical induction map associated with $\pi$}
to be the composite
\[
\can^\pi
:=
(\mk^A)^{-1}\circ\die^\pi.
\]
Thus
\[
\can^\pi\colon
{}\lRes{R\e}_{R\n}(M)
\longrightarrow
{}\lRes{R\e}_{R\n}(A_+)
\]
is an $R\e$-morphism.

Equivalently, after suppressing the restriction symbols, we have the
commutative triangle
\[
\xymatrix@R=4ex@C=4ex{
	A_+ \ar[rr]^{\mk^A} && A^+\\
	& M \ar[lu]^{\can^\pi}\ar[ru]_{\die^\pi} &
}
\]
of $R\e$-functors.
\end{Cons}
We now characterize when the canonical induction map
is a morphism of Green functors.
\begin{Prop}[Multiplicativity of canonical induction]
\label{prop:canonical-induction-multiplicative}
Let $M$ be a Green $R\n$-functor and let $A$ be a Green
$R\e$-functor. Suppose that
$\mk^A\colon A_+\isoto A^+$
is an isomorphism of Green $R\n$-functors, and let
$\pi\colon
{}\lRes{R\s}_{R\n}(M)
\to
{}\lRes{R\s}_{R\e}(A)$
be an $R\s$-morphism.
Then the following conditions are equivalent:
\begin{enumerate}
	\item $\pi$ is a Green $R\s$-morphism;
	\item $\die^\pi$ is a Green $R\e$-morphism;
	\item $\can^\pi$ is a Green $R\e$-morphism.
\end{enumerate}
\end{Prop}
\begin{proof}
The equivalence of the first two conditions follows from
\Cref{cor:upper-plus-Green-adjunction}, since the ordinary upper-plus
adjunction of \Cref{thm:upper-plus-adjunction} restricts to Green
functors.

Since $\mk^A$ is an isomorphism of Green $R\n$-functors, both
$\mk^A$ and $(\mk^A)^{-1}$ restrict to morphisms of Green
$R\e$-functors by \Cref{rem:Green-restrictions}. Hence
$\die^\pi$ is a Green $R\e$-morphism if and only if
$\can^\pi=(\mk^A)^{-1}\circ\die^\pi$ is one.
\end{proof}
For the explicit formula and the splitting criterion below, we return
to the hypothesis that $|G|$ is invertible in $R$ for every
$G\in\mathfrak K$.
The inverse mark formula gives an explicit description of the
canonical induction map.
\begin{Prop}[Explicit canonical induction formula]
\label{prop:canonical-induction-formula}
Assume that $|G|$ is invertible in $R$ for every
$G\in\mathfrak K$.
Let $M$, $A$, and $\pi$ be as in
\Cref{con:canonical-induction-map}.
For $G\in\mathfrak K$ and $m\in M(G)$, one has
\[
\can_G^\pi(m)
=
\frac{1}{|G|}
\sum_{U\leq V\leq G}
|U|\mu(U,V)
\left[
U,\,
{}\lres{U}_V
\left(
\pi_V({}\lres{V}_G(m))
\right)
\right]_G.
\]
Here the sum runs over all pairs of subgroups $U\leq V$ of~$G$.
\end{Prop}
\begin{proof}
Set $\xi=\die_G^\pi(m)$. By
\Cref{con:tom-Dieck-map},
\[
\xi^V=\pi_V\bigl({}\lres{V}_G(m)\bigr)
\qquad (V\leq G).
\]
Since $\can_G^\pi(m)=(\mk_G^A)^{-1}(\xi)$, the result follows
immediately from \Cref{prop:inverse-mark-formula}.
\end{proof}
In Boltje's classical setting, this specializes to
\cite[Eq.~(6.1a)]{Boltje98a}.
In classical canonical induction formulas, the projections
onto linear characters or exprojective $p$-permutation
modules need not preserve multiplication; see
\cite{Boltje98a,BarkerMutlu19}.
Consequently, the corresponding canonical induction
maps need not be morphisms of Green functors.

We now determine when the map $\can^\pi$ is a
canonical induction formula for $\lin^\nu$.
The following criterion generalizes
\cite[Proposition~6.4]{Boltje98a}.
\begin{Thm}[Canonical induction criterion]
\label{thm:canonical-induction-criterion}
Assume that $|G|$ is invertible in $R$ for every
$G\in\mathfrak K$.
Let $M$ be an $R\n$-functor and let $A$ be an
$R\e$-functor.
Suppose that we are given an $R\e$-morphism
$\nu\colon A\to{}\lRes{R\e}_{R\n}(M)$
and an $R\s$-morphism
$\pi\colon{}\lRes{R\s}_{R\n}(M)\to
{}\lRes{R\s}_{R\e}(A)$.
Then the following conditions are equivalent:
\begin{enumerate}
\item $\can^\pi$ is a canonical induction formula
for $\lin^\nu$.
\item For every $G\in\mathfrak K$ and
$m\in M(G)$,
\[
\nu_G(\pi_G(m))-m
\in
\sum_{H<G}\lind{G}_H(M(H)).
\]
\end{enumerate}
\end{Thm}
\begin{proof}
The argument is the same as that of
\cite[Prop.~6.4]{Boltje98a}; we indicate the points needed in the
present setting.

By
\Cref{prop:canonical-induction-formula,con:linearization-map} and
\Cref{lem:Burnside-action-compatibility},
\[
e_G^G\cdot
\bigl(
\lin_G^\nu\can_G^\pi(m)
\bigr)
=
e_G^G\cdot\nu_G(\pi_G(m))
\]
for every $G\in\mathfrak K$ and $m\in M(G)$, since $e_G^G$
annihilates every term induced from a proper subgroup.

If $\lin^\nu\circ\can^\pi=\id_M$, the displayed identity and
\Cref{prop:Burnside-decomposition} give
\[
\nu_G(\pi_G(m))-m\in\cI(M)(G).
\]
Conversely, assume~\textup{(b)}. Since $\lin^\nu$ is an
$R\n$-morphism and $\can^\pi$ is an $R\e$-morphism, the composite
$\lin^\nu\circ\can^\pi$ is an $R\e$-endomorphism and hence commutes
with restriction to subgroups. The induction argument of
\cite[Prop.~6.4]{Boltje98a}, using
\Cref{cor:Burnside-zero-detection}, therefore gives
$\lin_G^\nu\can_G^\pi=\id_{M(G)}$
for every $G\in\mathfrak K$.
\end{proof}
\begin{Rem}[The $p$-permutation formula]
\label{rem:p-permutation-canonical-induction}
The canonical induction formula of
\cite{BarkerMutlu19} is recovered from
\Cref{prop:canonical-induction-formula} and
	\Cref{thm:canonical-induction-criterion}
in the classical Mackey setting, where
$\mathcal R$ consists of injective homomorphisms.
In this case, $\mathcal S$ contains only isogations.

Take
\[
M(G)=\mathbb Q T(G),
\qquad
A(G)=\mathbb Q T^{\mathrm{ex}}(G),
\]
where $T(G)$ is the $p$-permutation ring and
$T^{\mathrm{ex}}(G)$ is its subring spanned by the
classes of exprojective modules.
The inclusion
$\nu\colon A\longrightarrow M$
is a morphism of Green restriction functors, while
$\pi_G$ is the projection onto
$\mathbb Q T^{\mathrm{ex}}(G)$ along the
indecomposable non-exprojective classes.
Since exprojectivity is preserved under isogation,
the projections $\pi_G$ commute with isogation maps.
Thus $\pi=(\pi_G)_{G\in\mathfrak K}$ is an
$R\s$-morphism.

Then \Cref{prop:canonical-induction-formula} gives
the map of \cite[Theorem~3.1]{BarkerMutlu19},
while \Cref{thm:canonical-induction-criterion}
reduces its splitting property to the
proper-induction condition verified in the proof
of that theorem.
\end{Rem}

\end{document}